\documentclass{article}
\usepackage[utf8]{inputenc}
\usepackage[
backend=biber,
style=alphabetic,
sorting=ynt
]{biblatex}
\usepackage{biblatex}
\usepackage{amsthm}
\usepackage{amsmath}
\usepackage{amssymb}
\usepackage{esint}
\usepackage{thmtools}
\usepackage{thm-restate}
\usepackage{graphicx}
\graphicspath{ {./images/} }

\usepackage{bbm}
\usepackage{color}
\usepackage[normalem]{ulem}
\newtheorem{theorem}{Theorem}[section]
\newtheorem{corollary}{Corollary}[theorem]
\newtheorem{lemma}[theorem]{Lemma}
\newtheorem{claim}[theorem]{Claim}

\newtheorem{definition}{Definition}[section]
\theoremstyle{remark}
\newtheorem*{remark}{Remark}
\DeclareMathOperator{\osc}{osc}
\DeclareMathOperator{\imm}{Imm}
\DeclareMathOperator{\supp}{supp}
\DeclareMathOperator{\identity}{Id}
\DeclareMathOperator{\eosc}{ess osc}
\DeclareMathOperator{\eim}{essR}
\DeclareMathOperator{\sign}{sgn}
\DeclareMathOperator{\trace}{Tr}
\DeclareMathOperator{\diam}{diam}
\DeclareMathOperator{\adj}{adj}
\DeclareMathOperator{\dist}{dist}
\DeclareMathOperator{\bdiv}{div}
\DeclareMathOperator{\convex}{Conv}

\title{Regularity of Finite Energy Deformations}
\author{Master thesis --- Daniel Rosenblatt\\ Advisor: Cy Maor}
\date{February 2023}
\begin{document}
\maketitle
\begin{abstract}
    We review results of {\v{S}}ver\'{a}k, and Goldstein-Haj{\l}asz-Pakzad on how to show the continuity of functions in a critical Sobolev space with positive Jacobian. In the final chapter we expand on the theory of $VMO$ functions, showing a version of the change of variables theorem in this regularity, which generalizes the first step in {\v{S}}ver\'{a}k's proof to a co-dimension 1 case.
\end{abstract}
\section{Introduction}
\subsection{Introduction}
Consider a body $\Omega\subset\mathbb{R}^n$ (of positive measure), and consider the deformation given by sending a particle located at position $x$ to $f(x)$ for some function $f:\Omega\rightarrow\mathbb{R}^n$, then we should ask what setting can we impose on the function $f$ such that the deformations can arise from physical situations. A first setting could be that of $W^{1,n}\left(\Omega;\mathbb{R}^n\right)$ so that we could give some measure of how much strain the body has under each deformation (we need at least one derivative to measure stretching). Unfortunately, under this single assumption we can still see non-physical deformations. One family of such deformations are those that have negative Jacobian in some set of positive measure, which would mean that in a sense, the body is being pulled inside-out, thus, a better setting would be the following:
$$\left\{f\in W^{1,n}\left(\Omega;\mathbb{R}^n\right)\mid \det\nabla f>0 \text{ a.e.\ }\right\}$$
These sort of settings comes from looking at finite energy functions of some elastic potential. A first example of such potential would be:
$$E\left(f\right)=\int\limits_\Omega \dist^n\left(\nabla f, SO\left(n\right)\right),$$
which clearly gives the first setting of $W^{1,n}\left(\Omega;\mathbb{R}^n\right)$. If we want the more physical setting as mentioned above, where the body cannot be pulled inside-out, we look at the potential:
$$E\left(f\right)=\int\limits_\Omega \dist^n\left(\nabla f, SO\left(n\right)\right)+\Theta\left(\det\nabla f\right)$$
Where $\Theta$ is some function that is infinite on all non-positive values, for example it can be chosen that $\Theta\left(t\right)=\log\left(t\right)$, which gives the setting of $W^{1,n}\left(\Omega;\mathbb{R}^n\right)$ and $\det\nabla f>0$ a.e.\\
An interesting question to ask about this setting is whether cavitation can occur, that is, whether such function can be discontinuous. Such question is natural both physically, in the sense of cavitation, and mathematically as the case $W^{1,n}$ is a borderline case of the embedding Theorems (the Sobolev Embedding theorem says that for every $\epsilon>0$ there is a continuous embedding $W^{1,n+\epsilon}\left(\Omega\right)\subset C\left(\Omega\right)$). Under additional assumptions (for example, in the statement below, we use extra regularity of the trace), we can show that cavitation does not occur, as a consequence of the following statement: \\
If $f\in W^{1,n}\left(\Omega;\mathbb{R}^n\right)$ with $\det\nabla f>0$ a.e.\ , $\Gamma=\partial\Omega$ is smooth, and also $\trace_\Gamma f\in W^{1,n}\left(\Gamma;\mathbb{R}^n\right)$ then 
$$\int\limits_\Omega g\left(f\left(x\right)\right) \det \nabla f\left(x\right)dx=\deg \left(f,\Omega ,y\right)\int\limits_{\mathbb{R}^n}g\left(z\right)dz$$
for any $g$ which is continuous and bounded.\\
Through this statement we can derive that functions in $W^{1,n}\left(\Omega;\mathbb{R}^n\right)$ with $\det\nabla f>0$ a.e.\ always have a continuous representative (which has a given form). The proof of this statement will be given in Chapter 2 (following in \cite{sverak}), and in Chapter 3 an alternative proof of this regularity result is given (following in \cite{goldstein}).\\
A natural step to follow, is to ask whether this results can be generalized to a case of co-dimension 1. To do this, the assumption of positive Jacobian needs to be modified to a relevant condition that can hold in co-dimension 1 (here the Jacobian is not defined).\\
Consider a sheet $\Omega\subset\mathbb{R}^n$, and take a deformation function $f:\Omega\rightarrow\mathbb{R}^{n+1}$ , then the equivalent setting to that of positive Jacobian is:
$$\imm\left(\Omega;\mathbb{R}^{n+1}\right)=\left\{f\in W^{1,n}\left(\Omega;\mathbb{R}^{n+1}\right)\mid \nu_f\text{ exists a.e.\ and }\nu_f\in W^{1,n}\left(\Omega;S^{n}\right) \right\}$$
Where $\nu_f\left(x\right)$ is the unit normal to $f\left(\Omega\right)$ at $f\left(x\right)$, and $S^n$ is the unit $n$-dimensional sphere. The similarity between the conditions can be seen by considering the function
$$\Tilde{F}\left(x,t\right)=f(x)+t\nu_f\left(x\right):\Omega\times\left[-d,d\right]\rightarrow\mathbb{R}^{n+1}$$
which under some additional assumptions has positive Jacobian. Another similarity is seen when we assume that $\nabla\nu_f=0$ a.e.\ In this setting, we can choose a coordinate system that gives $f\in W^{1,n}\left(\Omega;\mathbb{R}^n\right)$ and $\det\nabla f>0$ a.e. \\
The physical meaning of the set $\imm\left(\Omega;\mathbb{R}^{n+1}\right)$, is that each function in it has a finite deformation energy: the condition that the function is in $W^{1,n}$ means that the function has finite stretching energy, and the condition that the function has a normal in $W^{1,n}$ means that the function has finite bending energy. An example of such deformation energy is that of:
$$E\left(f\right)=\int\limits_\Omega \dist^n\left(df, O\left(n,n+1\right)\right)+\left\lvert\nabla\nu_f\right\rvert^n$$
Here the distance term measures stretching, while the other term measures bending.\\ 
In Chapter 4 we analyze the theory of degree introduced in Chapter 2 applied to functions with $VMO$ (Vanishing Mean Oscillation) regularity (as in the articles \cite{brenir1}, and \cite{brenir2}), we show that $\Tilde{F}\left(x,t\right)$ as above is in $VMO$, and 
generalize the change of variable formula introduced in Chapter 2 to functions of the same class as $\Tilde{F}\left(x,t\right)$. This already gives that cavitation does not occur with functions of this class.\\
Thus, a natural way to attack the question of whether $\imm\left(\Omega;\mathbb{R}^{n+1}\right)\subset C\left(\Omega;\mathbb{R}^{n+1}\right)$, is to try to mimic the proof in chapter 2 for $VMO(\mathbb{R}^n;\mathbb{R}^n)$ functions and thus obtain that the function $\Tilde{F}(x,t)$ as above is continuous, giving the desired result as a consequence. However, this is beyond the scope of this thesis.
\subsection{Notations}
\label{notation}
Before starting with the main argument of the paper, we introduce some notations preferences. Throughout this paper, we write $\Omega\subset\mathbb{R}^n$ as a body, that is a bounded connected set of positive measure, and we additionally assume that $\Gamma=\partial\Omega$ is a smooth $(n-1)$-dimensional manifold unless mentioned otherwise. Whenever we mention a function space without specifying the target space it can be assumed that if the domain is a body (usually $\Omega$), then the target space is of same dimension as the domain, for example $C\left(\Omega\right)$ is intended to be $C\left(\Omega;\mathbb{R}^n\right)$, whereas, if the domain is a sheet (usually $\Gamma$), then the target space is one dimension higher than the domain, for example, if $\Omega\subset\mathbb{R}^n$ and $\Gamma=\partial\Omega$, the notation $W^{1,n}\left(\Gamma\right)$ is intended to be $W^{1,n}\left(\Gamma;\mathbb{R}^n\right)$.\\
For some open set $A\subset\Omega$ with Lipschitz boundary, we denote by $\trace_{\partial A}$ the continuous extension from $C\left(\Omega\right)\cap W^{1,n}\left(\Omega\right)$ to $W^{1,n}\left(\Omega\right)$ of the restriction operator $f\rightarrow\left.f\right\rvert_{\partial A}$ with respect to the $W^{1,n}\left(A\right)$ norm.\\
Additionally, we denote by $T$ some tubular neighborhood of $\Gamma$ which has the following properties:
\begin{itemize}
    \item The projection map $P:T\rightarrow\Gamma$ is uniquely defined, and at least Lipschitz.
    \item For every $x\in T\backslash\Omega$ we have $2P\left(x\right)-x\in\Omega$.
    \end{itemize}
\begin{remark}
    Because $\Gamma$ is a smooth $\left(n-1\right)$-dimensional manifold, such $T$ can always be defined. If $\Gamma$ is not smooth but of the form $\Gamma'\times\left[-d,d\right]$, where $\Gamma'$ is smooth, then, we can still define $T$ for $\Gamma$ with the above properties.
\end{remark}
\section{Proof by Degree Theory}
\subsection{Idea}
The first proof we consider is one based on degree theory, in short this proof looks at the function, and by help of the degree, it can define a set of possible values for the function at a point, which maintains some sort of continuity with respect to all possible choices of representatives. A main result used is that this set contains a single value if and only if the function is continuous at the point.
\subsection{Introduction to Degree}
The main concept for this proof is that of the topological (Brouwer) degree. To better understand this concept, it is important to recall the results of Sard's Theorem.
\begin{theorem}[Sard]\cite[Chapters 2, 3]{milnor} 
Given some $f \in C^1\left(\Omega\right)$, if we define $S_f=f\left(\left\{x\in \Omega \mid \det \nabla f\left(x\right)=0\right\}\right)$, we get that $\lambda\left(S_f\right)=0$.
\end{theorem}
\begin{corollary}
 For $f\in C^1\left(\Omega\right) \cap C\left(\overline{\Omega}\right)$ and some $y\in f\left(\overline{\Omega}\right)\backslash\left(f\left(\Gamma\right)\cup S_f\right)$, we have that $0<\left\lvert f^{-1}\left(y\right)\right\rvert<\aleph_0$.
\end{corollary}
\begin{proof}\cite[Chapter 5.4]{ciarlet}
$f$ is a local diffeomorphism at every $x\in f^{-1}\left(y\right) $ (Inverse Function Theorem, $x\notin\Gamma$, and $\det \nabla f\left(x\right) \neq 0$). Hence, for every $x\in f^{-1}\left(y\right)$ there is a $x\in V_x\subset \Omega $ open, such that $\left.f\right\rvert_{V_x} \rightarrow \mathbb{R}^n$ is a diffeomorphism onto its image $W_x=f\left(V_x\right)$. Thus, by injectivity, $f^{-1}\left(y\right)\cap \left(V_x\backslash \left\{x\right\}\right)=\emptyset$, giving that $\bigcup\limits_{x\in f^{-1}\left(y\right)}V_x$ is an irreducible (no subcover) open cover of $f^{-1}\left(y\right)$. On the other hand $f^{-1}\left({y}\right)$ is a compact set (bounded as $\Omega$ is, and closed as preimage under continuous function (here uses $f\in C\left(\overline{\Omega}\right)$)), thus $f^{-1}\left(y\right)$ must be a finite set (as $\bigcup\limits_{x\in f^{-1}\left(y\right)}V_x$ must be a finite cover)).
\end{proof}
\begin{definition}
For a given function $f \in C^1\left(\Omega\right) \cap C\left(\overline{\Omega}\right)$, and a point
$y \in f\left(\Omega\right)\backslash\left(f\left(\Gamma\right)\cup S_f\right)$, we define the topological degree of f at y to be:
$$ \deg \left(f, \Omega ,y\right)= \sum \limits_{z \in f^{-1}\left(y\right)} \sign \det \nabla f\left(z\right)$$ 
And for $y\in \mathbb{R} ^n \backslash f\left(\overline{\Omega}\right)$ define $\deg \left(f,\Omega,y\right)=0$.
\end{definition}
\begin{remark}
This is well defined by the above theorem.
\end{remark}
Now we want to try extend the definition of degree to a more general setting.
\begin{claim}\cite[Chapter 5.4]{ciarlet}
Given $y,y'\in f\left(\overline{\Omega}\right)\backslash\left(f\left(\Gamma\right)\cup S_f\right)$ such that they both lie in the same connected component of $\mathbb{R}^n \backslash f\left(\Gamma\right)$, we have that $\deg \left(f,\Omega,y\right)=\deg \left(f,\Omega ,y'\right)$.
\end{claim}
Using the above Claim we can extend the definition to all of $\mathbb{R}^n \backslash f\left(\Gamma\right)$, by defining as before on $\mathbb{R}^n \backslash\left(f\left(\Gamma\right)\cup S_f\right)$, and for every $y\in S_f \backslash f\left(\Gamma\right)$, take any $y'$ in the same connected component of $\mathbb{R}^n\backslash f\left(\Gamma\right)$ as $y$, and define \\$\deg \left(f,\Omega ,y\right)=\deg \left(f,\Omega ,y'\right)$.
\begin{theorem}[Properties of deg]\cite[Chapter 5.4]{ciarlet}
\label{prop of deg}
For $\Omega \subset \mathbb{R}^n$ bounded, $f\in C^1\left(\Omega\right) \cap C\left(\overline{\Omega}\right)$, and $y\notin f\left(\Gamma\right)$, we have that:
\hfill
\begin{enumerate}
    \item $y \notin f\left(\overline{\Omega}\right) \Rightarrow \deg \left(f,\Omega ,y\right)=0$ or $\deg \left(f,\Omega ,y\right)\neq0 \Rightarrow y\in f\left(\Omega\right)$. 
    \item Given $y'\notin f\left(\Gamma\right)$ that lies in the same connected component of $\mathbb{R}^n \backslash f\left(\Gamma\right)$ as $y$, then we have $\deg \left(f,\Omega,y\right)=\deg \left(f,\Omega ,y'\right)$.
    \item There is an $\epsilon >0$ such that if $g\in C^1\left(\Omega\right) \cap C\left(\overline{\Omega}\right)$ with $\left\lVert f-g \right\rVert_{C\left(\bar{\Omega}\right)} \leq \epsilon$, then $y\notin g\left(\Gamma\right)$ and $\deg \left(f,\Omega ,y\right)=\deg \left(g,\Omega ,y\right)$. 
    \item Given a homotopy $t\in [0,1]\rightarrow f_t\in C^1\left(\Omega\right) \cap C\left(\overline{\Omega}\right)$ such that for all $t$, $y\notin f_t\left(\Gamma\right) $, then $\deg \left(f_t,\Omega ,y\right)=\deg \left(f_0,\Omega ,y\right)$.
    \item If $\Omega$ is connected, and $f$ is injective, then either $\deg\left(f,\Omega ,y\right)=1 $ for every $y\in f\left(\Omega\right)$, or $\deg\left(f,\Omega ,y\right)=-1 $  for every $y\in f\left(\Omega\right)$. 
    \item If $\Omega'\subset\Omega$, and $y\notin f\left(\Omega\backslash\Omega'\right)$  then $\deg\left(f,\Omega,y\right)=\deg\left(f,\Omega',y\right).$
\end{enumerate}
\end{theorem}
\begin{corollary}
\label{tr inv}
 If $f,g\in C^1\left(\Omega\right) \cap C\left(\overline{\Omega}\right)$ such that $\trace_\Gamma f=\trace_\Gamma g$, and $y\notin f\left(\Gamma\right)$ then $\deg \left(f,\Omega ,y\right)=\deg\left(g,\Omega, y\right)$.
\end{corollary}
\begin{proof}
Use part (4) on the Theorem together with the homotopy 
$$f_t=\left(1-t\right)\cdot f+t\cdot g$$
\end{proof}
 Now we can extend the definition of degree to the space $C\left(\overline{\Omega}\right)$. \\
 Given some $f\in C\left(\overline{\Omega}\right)$, and $y\notin f\left(\Gamma\right)$, take a sequence $f_n\in C^1\left(\Omega\right)\cap C\left(\overline{\Omega}\right)$ such that $\left\lVert f_n-f\right\rVert_{C\left(\overline{\Omega}\right)}\rightarrow 0$, and for every $n$, $y\notin f_n\left(\Gamma\right)$, and define 
 $$\deg \left(f, \Omega ,y\right)=lim_{n\rightarrow\infty}\deg \left(f_n,\Omega ,y\right)$$
 \begin{remark}
 By part (3) on Theorem \ref{prop of deg} this is well defined.
 \end{remark}
\begin{claim}
Theorem \ref{prop of deg} and its Corollary hold if we substitute $C^1\left(\Omega\right)\cap C\left(\overline{\Omega}\right)$ with $C\left(\overline{\Omega}\right)$.
\end{claim}
\begin{remark}
This is proved directly in \cite{ciarlet}, there the extension to $C\left(\overline{\Omega}\right)$ is done before proving the properties in Claim \ref{prop of deg}.
\end{remark}
We can additionally use Corollary \ref{tr inv} to extend the definition of degree to functions $f:\overline{\Omega}\rightarrow\mathbb{R}^n $ that have $\trace_\Gamma f\in C\left(\Gamma\right)$, at some $y\notin f\left(\Gamma\right)$, by taking any $\bar{f}\in C\left(\overline{\Omega}\right)$ with $\trace_\Gamma\bar{f}=\trace_\Gamma f$, and defining $\deg \left(f,\Omega ,y\right)=\deg \left(\bar{f},\Omega ,y\right)$.  
\begin{remark}
This last extension allows us to speak about the degree of continuous functions defined only at the boundary (by taking any extension $\bar{f}$ of it into the whole domain), thus we may use the following notation, at some $y\notin f\left(\Gamma\right)$ we write 
$$\deg \left(f,\Gamma, y\right)=\deg \left(\bar{f},\Gamma ,y\right)=\deg \left(\bar{f},\Omega ,y\right)=\deg \left(f,\Omega ,y\right).$$
\end{remark}
\begin{remark}
By Morrey's inequality, we have that if $f\in W^{1,n}\left(\Gamma\right)$ then $f \in C\left(\Gamma\right)$, hence, it has a defined degree (this is up to choosing the correspondent continuous representative of the trace, but whenever we use this notation we assume that the correct representative was chosen (for justification of this choice see Claim \ref{rep})).
\end{remark}
\begin{remark}
When we say that $f\in W^{1,n}\left(\Gamma\right)$ we mean it in the sense of Trace, meaning that $\trace_\Gamma f\in W^{1,n}\left(\Gamma\right)$, from here onward we will "abuse notation" and write it as above.
\end{remark}
We will now present a result, which will allow computing the more generalized definitions of degree.
\begin{theorem}[Change of Variable Formula]\cite[Proposition 3.1.2]{brouwer}
\label{changeofvar}
Let $f\in C^1\left(\Omega\right)\cap C\left(\overline{\Omega}\right)$, and let $g\in C^\infty\left(\mathbb{R}^n;\mathbb{R}\right)$ with support in some connected component of $\mathbb{R}^n\backslash f\left(\Gamma\right)$ which contains $y$, then:
\begin{equation}
\label{change of var}
    \int\limits_\Omega g\left(f\left(x\right)\right) \det \nabla f\left(x\right)dx=\deg \left(f,\Omega ,y\right)\int\limits_{\mathbb{R}^n}g\left(z\right)dz
\end{equation}
\end{theorem}
\begin{theorem}
\label{stok}
Let $f\in C^1\left(\overline{\Omega}\right)$ and let $g\in C^\infty\left(\mathbb{R}^n;\mathbb{R}\right)$ with support in some connected component of $\mathbb{R}^n\backslash f\left(\Gamma\right)$ which contains $y$.\\
Now consider a $C^2$ vector field $v=\left(v^1,...,v^n\right)$ on $\mathbb{R}^n$ such that $\bdiv v=g$, and define an $n$-form $\alpha$, and an $\left(n-1\right)$-form $\beta$ on $\mathbb{R}^n$ by $\alpha =g\cdot dz^1\wedge ...\wedge dz^n$ and $\beta =\sum\limits_{i=1}^n \left(-1\right)^{i-1}v^i\cdot dz^1\wedge ...\wedge \widehat{dz^i}\wedge...\wedge dz^n$, thus $d\beta =\alpha$, and 
$$\int\limits_\Omega g\left(f\left(x\right)\right)\det \nabla f\left(x\right)dx=\int\limits_\Omega f^*\left(\alpha\right)=\int\limits_\Omega f^*\left(d\beta\right)=\int\limits_\Omega df^*\left(\beta\right)=\int\limits_\Gamma f^*\left(\beta\right)$$
\end{theorem}
\begin{proof}
The equalities stated follow directly from definitions, the fact that pull-backs commute with exterior derivatives, and the use of Stokes Theorem.
\end{proof}
\begin{corollary} [Pull-back Formula]
Combining the last two Theorems we obtain that for $f\in C^1\left(\Gamma\right)$ and $g,\beta$ as before we have:
\begin{equation}
\label{eq deg}
    \deg \left(f,\Omega ,y\right)\int\limits_{\mathbb{R}^n}g\left(z\right)dz=\int\limits_\Gamma f^*\left(\beta\right)
\end{equation}
\end{corollary}
\begin{remark}
To Apply the Theorems, need to extend $f$ from $\Gamma$ to $\Omega$, but the formula is independent of the choice of extension.
\end{remark}
\subsection{Generalization of Properties of Degree}
Now we will see how to generalize some results obtained on the previous section.
\begin{lemma}
\label{eq deg g}
The Formula \eqref{eq deg} holds for $f\in W^{1,n}\left(\Gamma\right)$.
\end{lemma}
\begin{remark}
In \cite{sverak} this Lemma is stated for $f\in W^{1,p}\left(\Gamma\right)\cap C\left(\Gamma\right)$ for $p\geq n-1$, but, note that if $f\in W^{1,p}\left(\Gamma\right)$ for $p>n-1$, then $f$ already has a representative in $C\left(\Gamma\right)$.
\end{remark}
\begin{proof}
This is done by taking a sequence of functions $f_k\in C^1\left(\Gamma\right)$ such that $\lim\limits_{k\rightarrow\infty} \left\lVert f-f_k \right\rVert_{W^{1,n}\left(\Gamma\right)}=0$, by Morrey's inequality, this also guaranties that $\lim\limits_{k\rightarrow\infty} \left\lVert f-f_k\right\rVert_{C\left(\Gamma\right)}=0$. Now we can use the fact that the left side of \eqref{eq deg} is continuous with respect $C\left(\Omega\right)$ norm, and the right side is continuous with respect to $W^{1,n}\left(\Gamma\right)$ norm.
\end{proof}
\begin{theorem}
\label{change}
The Formula \eqref{change of var} holds for $g,y$ as before, and $f\in W^{1,n}\left(\Omega\right)$ such that also $f\in W^{1,n}\left(\Gamma\right)$.
\end{theorem}
\begin{remark}
Here $\deg \left(f,\Omega ,y\right)$ refers to $\deg \left(\trace_\Gamma f,\Gamma ,y\right)$ by choosing a continuous representative of $\trace_\Gamma f$, also $f\left(\Gamma\right)$ refers to the image under that representative.
\end{remark}
\begin{proof}
By the assumption on the domain, we can define the function $\nu :\Gamma\rightarrow \mathbb{S}^1$, sending each $x$ on $\Gamma$ to the unit outward normal of $\Omega$. For some $\epsilon>0$, we can define the diffeomorphism $h$ from $\Gamma\times \left(-2\epsilon,2\epsilon\right)$ to the tubular neighborhood $A_{2\epsilon}=\left\{z\in \mathbb{R}^n \mid \dist \left(z,\Gamma\right)<2\epsilon\right\}$, by $h\left(x,t\right)=x+t\cdot \nu \left(x\right)$. We note that for each $h\left(x,t\right)=z\in A_{2\epsilon}$, $x$ is the unique point in $\Gamma$ such that \\$\left\lvert t\right\rvert=\dist\left(z,x\right)=\dist\left(z,\Gamma\right)$.\\
Now define $\Omega_\epsilon=\left\{z\in\mathbb{R}^n\mid \dist\left(z,\Omega\right)<\epsilon\right\}$, let $D_\epsilon=\Omega_\epsilon\backslash\Omega$, and define the projection map $\pi :\Omega_\epsilon\rightarrow\Omega$ by: 
$$ \pi\left(z\right) =     \left\{ \begin{array}{rcl}
         z & \mbox{if}
         & z\in \overline{\Omega} \\ x  & \mbox{if} & h\left(x,t\right)=z\in D_\epsilon
                \end{array}\right.$$
Define $\bar{f}$ as the continuous representative of $\trace_\Gamma f$, and define the map: 
$$\tilde{f}\left(z\right)=\left\{ \begin{array}{rcl}
        f\left(z\right) & \mbox{if}
         & z\in \overline{\Omega} \\ \bar{f}\left(\pi\left(z\right)\right)  & \mbox{if} & z\in D_\epsilon
                \end{array}\right.$$
By using local coordinates on $\Gamma$, the definition of $\tilde{f}$, and the fact that $\bar{f}\in W^{1,n}\left(\Gamma\right)$, we obtain that $\tilde{f}\in W^{1,n}\left(\Omega_\epsilon\right)$.\\
Now consider any smooth function $g:\mathbb{R}^n\rightarrow\mathbb{R}$ compactly supported on the same connected component $W$ of $\mathbb{R}^n\backslash\bar{f}\left(\Gamma\right)$ as $y$, and consider the vector field $v$ as in Theorem \ref{stok} such that its derivatives are bounded (for example choose $v^i\left(w\right)=c_n\int\limits_{\mathbb{R}^n}\frac{w^i-z^i}{\left\lvert w-z\right\rvert^n}g\left(z\right)dz$ where $c_n$ is some suitable constant).\\
At this point, we need the following lemma:
\begin{lemma}
\label{inter}
for every $\xi\in Lip\left(\Omega_\epsilon;\mathbb{R}\right)$ with $\left.\xi\right\rvert_{\partial\Omega_\epsilon}=0$ we have:
\begin{equation}
\label{lem eq}
    \int\limits_{\Omega_\epsilon} \xi\left(x\right)g\left(\tilde{f}\left(x\right)\right) \det \nabla \tilde{f}\left(x\right)dx=-\int\limits_{\Omega_\epsilon} \xi_{,j}\left(x\right)v^k\left(\tilde{f}\left(x\right)\right) \left(\adj \nabla \tilde{f}\left(x\right)\right)^j_kdx
\end{equation}
Where the symbol $\xi_{,j}$ is the partial derivative of $\xi$ with respect to the $j$ coordinate
\end{lemma}
\begin{proof}
We may first observe that for any smooth function $\eta=\left(\eta^1,...,\eta^n\right)$ with compact support on $\Omega_\epsilon$, we have that: 
\begin{equation}
\label{ker}
    \int\limits_{\Omega_\epsilon} \eta^k_{,j}\left(x\right)\left(\adj \nabla \tilde{f}\left(x\right)\right)^j_kdx=0
\end{equation}
Indeed we get that by considering,
$$\frac{\partial\eta^k}{\partial x^j}\left(\adj\nabla\tilde{f}\right)^j_k=\sum\limits_{k=1}^n\det\nabla\left(\tilde{f}^1,...,\tilde{f}^{k-1},\eta^k,\tilde{f}^{k+1},...,\tilde{f}^n\right)$$
we get that: 
\begin{equation*}
    \begin{split}
        \int\limits_{\Omega_\epsilon} \eta^k_{,j}\left(x\right)\left(\adj \nabla \tilde{f}\left(x\right)\right)^j_kdx & = \int\limits_{\Omega_\epsilon}\sum\limits_{k=1}^n\det\nabla\left(\tilde{f}^1,...,\tilde{f}^{k-1},\eta^k,\tilde{f}^{k+1},...,\tilde{f}^n\right)\left(x\right)dx\\
        & = \int\limits_{\Omega_\epsilon}\sum\limits_{k=1}^n d\tilde{f}^1\wedge...\wedge d\tilde{f}^{k-1}\wedge d\eta^k\wedge d\tilde{f}^{k+1}\wedge...\wedge d\tilde{f}^n
    \end{split}
\end{equation*}
and the last integral is always zero by Stokes Theorem as $\eta$ is compactly supported on $\Omega_\epsilon$.\\
Now using the fact that $\left(\adj\nabla\tilde{f}\right)^j_k\in L^{\frac{n}{n-1}}\left(\Omega_\epsilon\right)$ we can extend Formula \eqref{ker} to any $\eta\in W^{1,n}_0\left(\Omega_\epsilon\right)$ (by Hölder inequality). \\
We additionally note that $\eta=\xi\cdot\left(v\circ\tilde{f}\right)\in W^{1,n}_0\left(\Omega_\epsilon\right)$ thus by Formula \eqref{ker} we get that:
\begin{equation*}
\begin{split}
    & -\int\limits_{\Omega_\epsilon} \xi_{,j}\left(x\right)v^k\left(\tilde{f}\left(x\right)\right) \left(\adj \nabla \tilde{f}\left(x\right)\right)^j_kdx\\
    &\qquad =\int\limits_{\Omega_\epsilon} \xi\left(x\right)\frac{\partial}{\partial x^j}\left(v^k\circ\tilde{f}\right)\left(x\right) \left(\adj \nabla \tilde{f}\left(x\right)\right)^j_kdx\\
    &\qquad =\int\limits_{\Omega_\epsilon} \xi\left(x\right)\frac{\partial v^k}{\partial x^l}\frac{\partial\tilde{f}^l}{\partial x^j}\left(x\right) \left(\adj \nabla \tilde{f}\left(x\right)\right)^j_kdx=\int\limits_{\Omega_\epsilon} \xi\left(x\right)\frac{\partial v^k}{\partial x^l}\delta^l_k \det \nabla \tilde{f}\left(x\right)dx\\
    &\qquad =\int\limits_{\Omega_\epsilon} \xi\left(x\right)g\left(\tilde{f}\left(x\right)\right) \det \nabla \tilde{f}\left(x\right)dx
\end{split}
\end{equation*}
as needed. Note that above, the first equality is the use of \eqref{ker} ,and the third equality follows from the definition of the determinant.
\end{proof}
Proof of Theorem \ref{change} (continued):\\
Now, consider the lipshitz function $\xi\left(x\right)=1-\frac{1}{\epsilon}\dist\left(x,\Omega\right)$. We note that it satisfies the conditions for Lemma \ref{inter}. Looking at the left side of \eqref{lem eq} we get: 
\begin{equation}
\label{lhs}
\begin{split}
    \int\limits_{\Omega_\epsilon} \xi\left(x\right)g\left(\tilde{f}\left(x\right)\right) \det \nabla \tilde{f}\left(x\right)dx & =\int\limits_{\Omega} g\left(\tilde{f}\left(x\right)\right) \det \nabla \tilde{f}\left(x\right)dx\\
    & =\int\limits_\Omega g\left(f\left(x\right)\right) \det \nabla f\left(x\right)dx
\end{split}
\end{equation}
Because $\det\nabla \tilde{f}\left(x\right)$ is a.e. zero in $D_\epsilon$, and inside $\Omega$ we have $\xi=1$, and $\tilde{f}=f$.\\
Next, we will rewrite $-\int\limits_{\Omega_\epsilon} \xi_{,j}\left(x\right)v^k\left(\tilde{f}\left(x\right)\right) \left(\adj \nabla \tilde{f}\left(x\right)\right)^j_kdx$:\\
First we notice that that $\xi$ is constant on $\Omega$, thus $\xi_{,j}=0$ on $\Omega$, hence, we can change the domain of integration to $D_\epsilon$. By following the definitions we get: 
\begin{equation*}
    \begin{split}
        & d\xi\wedge\tilde{f}^*\left(\beta\right)\\
        &\qquad =d\xi\wedge\sum\limits_{k=1}^n\left(-1\right)^{k-1}v^k\left(\tilde{f}\right)d\left(x^1\circ\tilde{f}\right)\wedge...\wedge \widehat{d\left(x^k\circ\tilde{f}\right)}\wedge...\wedge d\left(x^n\circ\tilde{f}\right)\\
&\qquad =d\xi\wedge\sum\limits_{k=1}^n\left(-1\right)^{k-1}v^k\left(\tilde{f}\right)d\tilde{f}^1\wedge...\wedge \widehat{d\tilde{f}^k}\wedge...\wedge d\tilde{f}^n\\
&\qquad =\sum\limits_{k=1}^n\left(-1\right)^{k-1}v^k\left(\tilde{f}\right)\det\nabla\left(\xi,\tilde{f}^1,..., \widehat{\tilde{f}^k},...,\tilde{f}^n\right)dx^1\wedge...\wedge dx^n\\
&\qquad =\sum\limits_{k=1}^n\left(-1\right)^{k-1}v^k\left(\tilde{f}\right)\xi_{,j}\left(\adj\nabla\left(\xi,\tilde{f}^1,..., \widehat{\tilde{f}^k},...,\tilde{f}^n\right)\right)^j_1dx^1\wedge...\wedge dx^n\\
&\qquad =\sum\limits_{k=1}^nv^k\left(\tilde{f}\right)\xi_{,j}\left(\adj\nabla\left(\tilde{f}^1,..., \tilde{f}^{k-1},\xi,\tilde{f}^{k+1},...,\tilde{f}^n\right)\right)^j_kdx^1\wedge...\wedge dx^n\\
&\qquad =v^k\left(\tilde{f}\right)\xi_{,j}\left(\adj\nabla\left(\tilde{f}\right)\right)^j_kdx^1\wedge...\wedge dx^n.
    \end{split}
\end{equation*}
Thus we obtain that:
\begin{equation}
\label{rhs1}
    -\int\limits_{\Omega_\epsilon} \xi_{,j}\left(x\right)v^k\left(\tilde{f}\left(x\right)\right) \left(\adj \nabla \tilde{f}\left(x\right)\right)^j_kdx=-\int\limits_{D_\epsilon} d\xi\wedge\tilde{f}^*\left(\beta\right)
\end{equation}
As in Lemma \ref{eq deg g} consider a sequence $\bar{f}_r\in C^1\left(\Gamma\right)$ converging to $\bar{f}$ in $W^{1,n}\left(\Gamma\right)$ norm (and thus also uniformly by Morrey), and consider on $D_\epsilon$, the functions defined by $\tilde{f}_r=\bar{f}\circ\pi$, and we can see that $\tilde{f}_r^*\left(\beta\right)\rightarrow\tilde{f}^*\left(\beta\right)$ and $\bar{f}_r^*\left(\beta\right)\rightarrow\bar{f}^*\left(\beta\right)$ on the space of $\left(n-1\right)$-forms with $L^p$ components over $D_\epsilon$ and $\Gamma$ respectively (the norm is the sum of $L^p$ norms of the components of the form). Finally for $t\in \left(0,\epsilon\right)$ we define $\Gamma_t=\left\{x\in\mathbb{R}^n\mid\dist\left(x,\Omega\right)=t\right\}$, thus $D_\epsilon=\bigcup\limits_{t\in \left(0,\epsilon\right)}\Gamma_t$, by the Fubini-Tonelli Theorem, and the fact that on $D_\epsilon$, $d\xi=-\frac{1}{\epsilon}dt$, we obtain that: 
\begin{equation}
\label{rhs2}
    \begin{split}
    -\int\limits_{D_\epsilon} d\xi\wedge\tilde{f}^*\left(\beta\right) & =\lim\limits_{r\rightarrow\infty}-\int\limits_{D_\epsilon} d\xi\wedge\tilde{f}_r^*\left(\beta\right)=\lim\limits_{r\rightarrow\infty}\int\limits_0^\epsilon\frac{1}{\epsilon}\left(\int\limits_{\Gamma_t} \tilde{f}_r^*\left(\beta\right)\right)dt\\
    & =\lim\limits_{r\rightarrow\infty}\int\limits_0^\epsilon\frac{1}{\epsilon}\left(\int\limits_{\Gamma_t} \pi^*\left(\bar{f}_r^*\left(\beta\right)\right)\right)dt=\lim\limits_{r\rightarrow\infty}\int\limits_0^\epsilon\frac{1}{\epsilon}\left(\int\limits_{\Gamma} \bar{f}_r^*\left(\beta\right)\right)dt\\
    & =\lim\limits_{r\rightarrow\infty}\int\limits_{\Gamma} \bar{f}_r^*\left(\beta\right)=\int\limits_{\Gamma} \bar{f}^*\left(\beta\right)=\deg \left(f,\Omega ,y\right)\int\limits_{\mathbb{R}^n}g\left(z\right)dz
    \end{split}
\end{equation}
Where the last equality is exactly the result of Lemma \ref{eq deg g}.\\
Thus combining the result of Lemma \ref{inter} with Equations \eqref{lhs},\eqref{rhs1}, and \eqref{rhs2}, we obtained the desired equality:
$$\int\limits_\Omega g\left(f\left(x\right)\right) \det \nabla f\left(x\right)dx=\deg \left(f,\Omega ,y\right)\int\limits_{\mathbb{R}^n}g\left(z\right)dz$$
\end{proof}
\begin{corollary}\cite{sverak} (not used in main proof).
Let $f\in W^{1,n}\left(\Omega\right)$ such that $f\in W^{1,n}\left(\Gamma\right)$, and $g\in C\left(\Omega\right)$ bounded.
Then, if the measure of $f\left(\Gamma\right)$ is zero, we have:
$$\int\limits_\Omega g\left(f\left(x\right)\right) \det \nabla f\left(x\right)dx=\int\limits_{\mathbb{R}^n}\deg \left(f,\Omega ,y\right)g\left(y\right)dy$$
and in particular if $g\equiv1$
$$\int\limits_\Omega  \det \nabla f\left(x\right)dx=\int\limits_{\mathbb{R}^n}\deg \left(f,\Omega ,y\right)dy$$
\end{corollary}
To show this the following lemma is used:
\begin{lemma}\cite{sverak}
Let $f\in W^{1,1}\left(\Omega\right)$ and let $A=\left\{x\in\Omega\mid\det\nabla f\left(x\right)\neq0\right\}$\\
Now let some $M\subset\mathbb{R}^n$ of zero measure and $B\subset\Omega$ measurable such that for a.e.\ $x\in B$ we have $f\left(x\right)\in M$. Then, the measure of $A\cap B$ is zero.
\end{lemma}
\subsection{Useful Definitions Involving Degree}
Now we can consider the following set:
\begin{definition} For $f\in W^{1,n}\left(\Omega\right)$ such that $f\in W^{1,n}\left(\Gamma\right)$ define:
$$E\left(f,\Omega\right)=f\left(\Gamma\right)\cup\left\{y\in \mathbb{R}^n\backslash f\left(\Gamma\right)\mid\deg\left(f,\Omega,y\right)\geq1\right\}$$
Note that here, $f\left(\Gamma\right)$ is in the sense of picking a continuous representative of $\trace_\Gamma f$.
\end{definition}
\begin{remark}
Notice that $E\left(f,\Omega\right)$ is bounded as the degree is non-zero in a bounded set, and also closed as the degree is constant on connected components of $\mathbb{R}^n\backslash f\left(\Gamma\right)$, and $E\left(f,\Omega\right)$ contains $f\left(\Gamma\right)$. Thus $E\left(f,\Omega\right)$ is compact.
\end{remark}
\begin{remark}
Note that $E\left(f,\Omega\right)$ is also connected as the union of some connected components of $\mathbb{R}^n\backslash f\left(\Gamma\right)$ and $f\left(\Gamma\right)$.
\end{remark}
\begin{corollary}
\label{belong}
If $\det\nabla f>0$ a.e.\ in $\Omega$ then $f\left(x\right)\in E\left(f,\Omega\right)$ for a.e.\ $x\in\Omega$.
\end{corollary}
\begin{proof}
Let $x_0\in \Omega$, if $f\left(x_0\right)\in f\left(\Gamma\right)$ we are done. Else consider $W_{x_0}$ to be the connected component of $\mathbb{R}^n\backslash f\left(\Gamma\right)$ containing $f\left(x_0\right)$, and a sequence of smooth functions $g_i$ converging to $\mathbbm{1}_{W_{x_0}}$ monotonically and in $L^1$ norm, then we have:
\begin{equation*}
\begin{split}
    \deg\left(f,\Omega,f\left(x_0\right)\right)\lambda_n\left(W_{x_0}\right) & =\deg\left(f,\Omega,f\left(x_0\right)\right)\int\limits_{\mathbb{R}^n}\mathbbm{1}_{W_{x_0}}\left(z\right)dz\\
    & =\lim\limits_{i\rightarrow\infty}\deg\left(f,\Omega,f\left(x_0\right)\right)\int\limits_{\mathbb{R}^n}g_i\left(z\right)dz\\
    & =\lim\limits_{i\rightarrow\infty}\int\limits_\Omega g_i\left(f\left(x\right)\right) \det \nabla f\left(x\right)dx\\
    & =\int\limits_\Omega \mathbbm{1}_{W_{x_0}}\left(f\left(x\right)\right) \det \nabla f\left(x\right)dx,
\end{split}
\end{equation*}
where the last equality follows from dominated convergence. Now consider the set $U=\left\{x\in\Omega\backslash f^{-1}\left(f\left(\Gamma\right)\right)\mid \lambda_n\left(f^{-1}\left(W_x\right)\right)=0\right\}$, and using that $\mathbb{R}^n\backslash f\left(\Gamma\right)$ has at most a countable number of connected components, we obtain that:
\begin{equation*}
\begin{split}
    \lambda_n\left(U\right) & \leq\lambda_n\left(f^{-1}\left(f\left(U\right)\right)\right)=\lambda_n\left(f^{-1}\left(\bigcup\limits_{x\in U}f\left(x\right)\right)\right)\leq \lambda_n\left(f^{-1}\left(\bigcup\limits_{x\in U}W_x\right)\right)\\
    & =\lambda_n\left(\bigcup\limits_{x\in U}f^{-1}\left(W_x\right)\right)=\lambda_n\left(\bigcup\limits_{W s.t.\exists x \in U, W_x=W }f^{-1}\left(W\right)\right)=0,
\end{split}
\end{equation*}
where the last equality holds as a countable union (countable connected components of $\mathbb{R}^n\backslash f\left(\Gamma\right)$) of zero measure set. Thus U is of zero measure, implying that for a.e.\ $x\in\Omega\backslash f^{-1}\left(f\left(\Gamma\right)\right)$ we have that $\lambda_n\left(f^{-1}\left(W_x\right)\right)\neq0$,
hence, using $\det\nabla f>0$ a.e.\ in $\Omega$ we obtain that:
$$\deg\left(f,\Omega,y\right)\lambda_n\left(W_{x_0}\right)=\int\limits_\Omega \mathbbm{1}_{W_{x_0}}\left(f\left(x\right)\right) \det \nabla f\left(x\right)dx>0.$$ 
Using the fact that the degree is always an integer we indeed obtain that 
$\deg\left(f,\Omega,f\left(x\right)\right)\geq 1$ for a.e.\ $x\in\Omega\backslash f^{-1}\left(f\left(\Gamma\right)\right)$ thus $f\left(x\right)\in E\left(f,\Omega\right)$ for a.e.\ $x\in\Omega$ as needed.
\end{proof}
\begin{remark}
From now on we will assume that $\det\nabla f>0$ a.e.\ in $\Omega$.
\end{remark}
\begin{lemma}[Monotonicity of $E\left(f,\cdot\right)$]
\label{mono E}
Let $\Omega_1\subset\Omega_2\subset\overline{\Omega_2}\subset\Omega$ all with smooth boundaries, if $\Gamma_i=\partial\Omega_i$, and $f\in W^{1,n}\left(\Gamma_i\right)$ then $E\left(f,\Omega_1\right)\subset E\left(f,\Omega_2\right)$.
\end{lemma}
\begin{proof}
For every  $x_0\in \Gamma_1$ consider a small segment $\gamma_{x_0}$ through $x_0$ such that $\gamma_{x_0}\subset{\overline{\Omega_2}}$, and for each pair $x_0,x_1\in\Gamma_1$, the segments $\gamma_{x_0},\gamma_{x_1}$ are disjoint (here we use that $\Gamma_1$ is smooth and thus has a tubular neighborhood). If we have that $f$ is absolutely continuous on $\gamma_{x_0}$, and $f\left(x_0\right)\notin E\left(f,\Omega_2\right)$, then, using that $E\left(f,\Omega_2\right)$ is closed, we have that $f\left(\gamma_{x_0}\right)\cap E\left(f,\Omega_2\right)=\emptyset$ (after restricting $\gamma_{x_0}$ to some smaller segment of itself). Thus, if we assume by contradiction that $\mathcal{H}_{n-1}\left(\left\{x\in\Gamma_1\mid f\left(x\right)\notin E\left(f,\Omega_2\right)\right\}\right)>0$, and use the fact that for a.e.\ $x\in \Gamma_1$ (w.r.t.\ the $\left(n-1\right)$-dimensional Hausdorff measure on $\Gamma_1$) we have that $f$ is absolutely continuous on $\gamma_x$ (proof identical to that in \cite[Theorem 4.21]{evans}), then, by using the Fubini-Tonelli Theorem on the indicator function, we obtain that $\lambda_{n}\left(\left\{x\in\Omega_2\mid f\left(x\right)\notin E\left(f,\Omega_2\right)\right\}\right)>0$, which is a contradiction to Corollary \ref{belong}.\\
Thus $f\left(x\right)\in E\left(f,\Omega_2\right)$ for a.e.\ $x\in\Gamma_1$ (with respect to surface measure), and by choosing the appropriate continuous representative and using compactness of $E\left(f,\Omega_2\right)$, we have that $f\left(\Gamma_1\right)\subset E\left(f,\Omega_2\right)$.\\
Now let $y\in E\left(f,\Omega_1\right)$. We just saw that if $y\in f\left(\Gamma_1\right)$, then $y\in E\left(f,\Omega_2\right)$. Also clearly if $y\in f\left(\Gamma_2\right)$, then $y\in E\left(f,\Omega_2\right)$.\\
Else take some smooth function $g$ supported in the same connected components of $\mathbb{R}^n\backslash f\left(\Gamma_1\right)$ and $\mathbb{R}^n\backslash f\left(\Gamma_2\right)$ as $y$ with $\int\limits_{\mathbb{R}^n} g\left(z\right)dz=1$, then we have that:
\begin{equation*}
    \begin{split}
        \deg\left(f,\Omega_2,y\right) & =\deg\left(f,\Omega_2,y\right)\int\limits_{\mathbb{R}^n} g\left(z\right)dz=\int\limits_{\Omega_2} g\left(f\left(x\right)\right) \det \nabla f\left(x\right)dx\\
& \geq\int\limits_{\Omega_1} g\left(f\left(x\right)\right) \det \nabla f\left(x\right)dx=\deg\left(f,\Omega_1,y\right)\int\limits_{\mathbb{R}^n} g\left(z\right)dz \\
& =\deg\left(f,\Omega_1,y\right)\geq 1
    \end{split}
\end{equation*}
Where the last inequality follows from $y\in E\left(f,\Omega_1\right)$.\\
Thus we see that $y\in E\left(f,\Omega_2\right)$ as needed.
\end{proof}
Now let $r_a$ be some positive number s.t.\ $B_{r_a}\left(a\right)\subset\Omega$. With this we can state the following Claims:
\begin{claim}
\cite[Chapter 4.3]{evans}
\label{rep tr}
Let $a\in\Omega$, and $r\in \left(0,r_a\right)$, then we have that the function defined for any $x\in S_r\left(a\right)$ by:
$$\hat{f}_r\left(x\right) = \lim\limits_{\underset{\rho\in(0,r_a)}{\rho\rightarrow0}}\frac{1}{\lambda_n\left(B_\rho\left(x\right)\cap B_r\left(a\right)\right)}\int\limits_{B_\rho\left(x\right)\cap B_r\left(a\right)}f(z)dz,$$
is a representative of $\trace_{S_r\left(a\right)}$.
\end{claim}
\begin{claim}
\label{rep}
Let $a\in\Omega$, and for any $r\in \left(0,r_a\right)$ consider $\hat{f}_{a,r}$ to be the representative, as in Claim \ref{rep tr}, of $\trace_{S_r\left(a\right)} f$, and for any $x\in S_r\left(a\right)$ define $\hat{f_a}\left(x\right)=\hat{f}_{a,r}\left(x\right)$. Then, $\hat{f_a}$ is a representative of $f$ on $B_{r_a}\left(a\right).$ 
\end{claim}
To prove this, we first state the following definition:
\begin{definition}[Nicely Shrinking Family]
The family of sets $\left\{E_i\right\}$ is said to be nicely shrinking to $x$ if there is a sequence of balls $B_{r_i}\left(x\right)$ satisfying that $\inf r_i=0$, $E_i\subset B_{r_i}\left(x\right)$ and $\sup \frac{\lambda_n\left(B_{r_i}\left(x\right)\right)}{\lambda\left(E_i\right)}<\infty$.
\end{definition}
We recall the Theorem:
\begin{theorem}
\cite[Theorem 7.10]{rudin}
    Let $f\in L^1\left(\Omega\right)$, if we associate to each $x\in\Omega$ a family of Borel sets $\left\{E_i\left(x\right)\right\}$ that nicely shrinks to $x$. Then
    $$f\left(x\right)=\lim\limits_{E_i\rightarrow x}\frac{1}{\lambda_n\left(E_i\left(x\right)\right)}\int\limits_{E_i\left(x\right)}fd\lambda_n$$
    at every Lebesgue point of f, and in particular, for a.e.\ x.
\end{theorem}
\begin{remark}
In \cite{rudin} they do this with sequences, and not families, but the proof is the same.
\end{remark}
\begin{proof} (of Claim \ref{rep}).
All that remains in this proof, is to see that for each $x\in B_{r_a}\left(a\right)$, the family $\left\{E_\rho\left(x\right)=B_\rho\left(x\right)\cap B_{R}\left(a\right)\right\}_\rho$, where $R=\dist\left(x,a\right)$, nicely shrinks to $x$.\\
We can note that if $\rho\geq2R$, then $E_\rho\left(x\right)=B_R\left(a\right)\subset B_{2R}\left(x\right)$, and thus take $r_\rho=2R$, and get $$\frac{\lambda_n\left(B_{r_\rho}\left(x\right)\right)}{\lambda\left(E_\rho\right)}=\frac{\lambda_n\left(B_{2R}\left(x\right)\right)}{\lambda\left(B_R\left(a\right)\right)}=2^n.$$
If $\rho<2R$, we know that $E_\rho\left(x\right)\subset B_\rho\left(x\right)$, thus we take $r_\rho=\rho$, which gives $\lambda_n\left(B_\rho\left(x\right)\right)=C\rho^n$, where $C$ is some constant depending on the dimension $n$.\\
If $\rho<\sqrt{2}R$, we can notice that (see Figure \ref{fig:my_label} bellow) $S_\rho\left(x\right)\cap S_R\left(a\right)$ is an $\left(n-1\right)$-dimensional ball, on the hyperplane perpendicular to the line containing both $a$ and $x$, with a radius given by: $$r'_\rho=\rho\sqrt{1-\left(\frac{\rho}{2R}\right)^2}$$
We observe that the $n$-dimensional cones with base as the above mentioned ball, and vertices being $x$ and the intersection point of the ray from $x$ towards $a$ and $S_\rho\left(x\right)$ are both contained in $E_\rho\left(x\right)$, thus we have the bound:
$$\lambda_n\left( E_\rho\left(x\right)\right)\geq C'\rho r_\rho'^{n-1}$$ 
Where $C'$ is a constant depending on dimension $n$, giving the bound:
\begin{equation*}
    \begin{split}
        \frac{\lambda_n\left(B_{r_\rho}\left(x\right)\right)}{\lambda\left(E_\rho\right)} & =\frac{\lambda_n\left(B_{\rho}\left(x\right)\right)}{\lambda\left(E_\rho\right)}\leq \frac{C\rho^n}{C'\rho r_\rho'^{n-1}}=\frac{C}{C'}\left(\frac{\rho}{r_\rho'}\right)^{n-1}\\
        & =\frac{C}{C'}\left(\frac{1}{1-\left(\frac{\rho}{2R}\right)^2}\right)^{\frac{n-1}{2}}<\frac{C}{C'}\left(\frac{1}{1-\left(\frac{\sqrt{2}R}{2R}\right)^2}\right)^{\frac{n-1}{2}}=\frac{C}{C'}2^{\frac{n-1}{2}}
    \end{split}
\end{equation*}
If $2R>\rho\geq\sqrt{2}R$ we notice that (see Figure \ref{fig:my_label} bellow) the $n$-dimensional cones , with base as the $\left(n-1\right)$-dimensional ball of radius $R$ contained in the hyperplane passing through $a$, and perpendicular to the line through $x$ and $a$, and heights $R$ (towards $x$), and $\rho-R$ (in the direction opposite of $x$), are both contained in $E_\rho\left(x\right)$, thus we have the bound:
$$\lambda_n\left( E_\rho\left(x\right)\right)\geq C'\rho R^{n-1}$$
Giving the bound:
$$\frac{\lambda_n\left(B_{r_\rho}\left(x\right)\right)}{\lambda\left(E_\rho\right)}=\frac{\lambda_n\left(B_{\rho}\left(x\right)\right)}{\lambda\left(E_\rho\right)}\leq \frac{C\rho^n}{C'\rho R^{n-1}}=\frac{C}{C'}\left(\frac{\rho}{R}\right)^{n-1}<\frac{C}{C'}2^{n-1}$$
Thus we get the global bound: $$\sup \frac{\lambda_n\left(B_{r_\rho}\left(x\right)\right)}{\lambda\left(E_\rho\right)}<\max\left\{\frac{C}{C'}2^{n-1},\frac{C}{C'}2^{\frac{n-1}{2}},2^n\right\}<\infty$$
As required.
\begin{figure}[h]
    \centering
    \includegraphics[width=\textwidth]{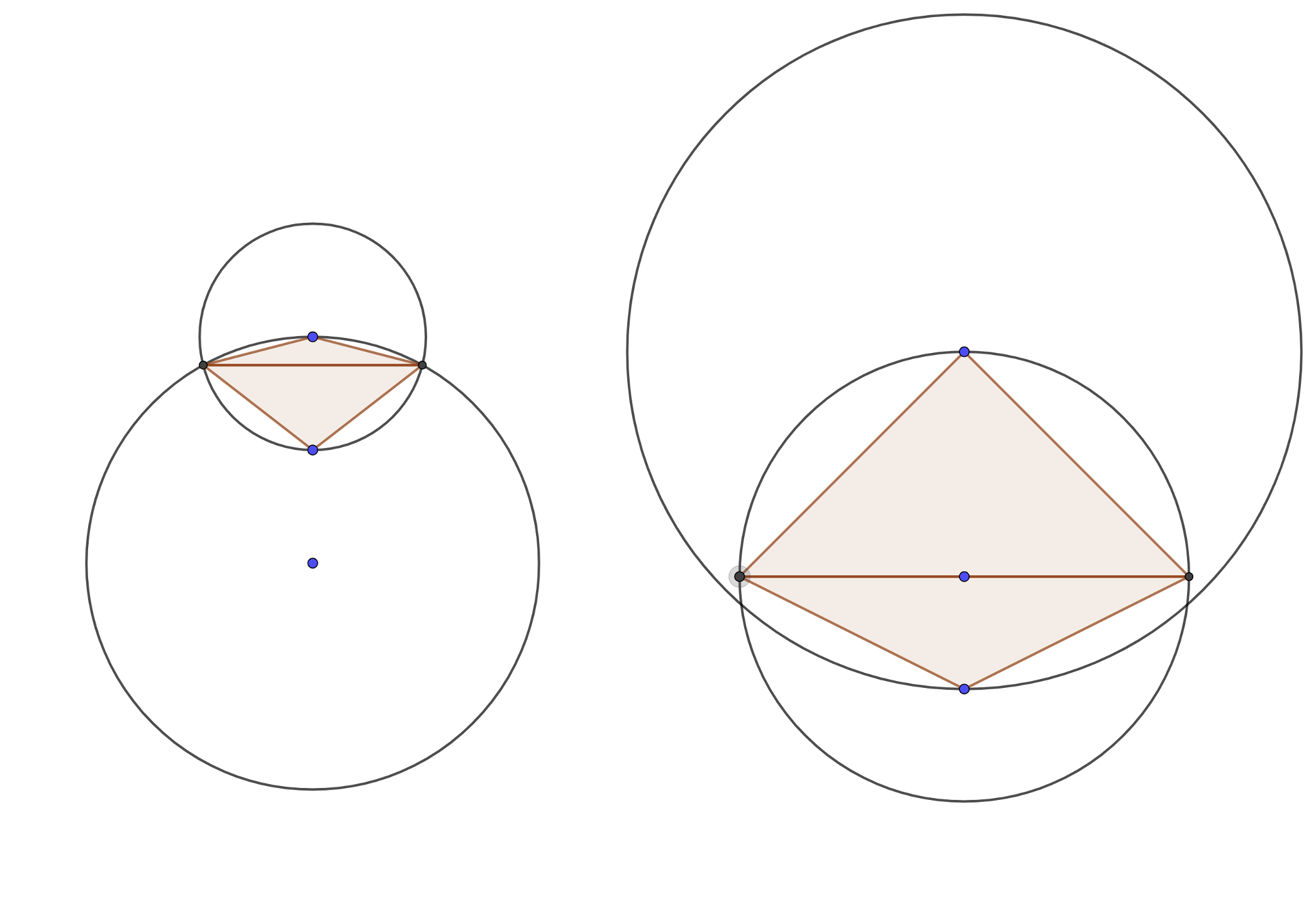}
    \caption{Left: the case $\rho<\sqrt{2}R$ Right: the case $\sqrt{2}R<\rho<2R$}
    \label{fig:my_label}
\end{figure}
\end{proof}
\begin{corollary}
\label{cor rep}
Let $a\in\Omega$, and $\bar{f}$ be any representative of $f$ then the set $Z_{a,\bar{f}}=\left\{r\in\left(0,r_a\right)\mid\left.\bar{f}\right\rvert_{S_r\left(a\right)} \text{ is not a representative of }\trace_{S_r\left(a\right)}f\right\}$ is of zero measure.
\end{corollary}
\begin{proof}
Assume by contradiction that the set $Z_{a,\bar{f}}$ has positive measure, then using the Minkowski inequality and the Fubini-Tonelli Theorem we see the following:
\begin{equation*}
    \begin{split}
        0 & <\int\limits_{Z_{a,\bar{f}}}\left(\int\limits_{S_r\left(a\right)}\left\lvert \trace_{S_r\left(a\right)}f-\bar{f}\right\rvert d\sigma_{S_r\left(a\right)}\right)dr\\
        & \leq\int\limits_0^{r_a}\left(\int\limits_{S_r\left(a\right)}\left\lvert \trace_{S_r\left(a\right)}f-\bar{f}\right\rvert d\sigma_{S_r\left(a\right)}\right)dr=\int\limits_0^{r_a}\left(\int\limits_{S_r\left(a\right)}\left\lvert\hat{f}-\bar{f}\right\rvert d\sigma_{S_r\left(a\right)}\right)dr\\
        & =\int\limits_{B_{r_a}\left(a\right)}\left\lvert\hat{f}-\bar{f}\right\rvert dz\leq\int\limits_{B_{r_a}\left(a\right)}\left\lvert\hat{f}-f\right\rvert dz+\int\limits_{B_{r_a}\left(a\right)}\left\lvert f-\bar{f}\right\rvert dz=0+0=0
    \end{split}
\end{equation*}
Giving a clear contradiction.
\end{proof}
For the next definition we need the following Claim:
\begin{claim}
\label{ineq of surf}
Let $a\in \Omega$, let $f\in W^{1,n}\left(\Omega\right)$, and $r\in (0,r_a)$, then for a.e.\ $\rho\in \left(0,r\right)$ we have that $f\in W^{1,n}\left(S_\rho \left(a\right)\right)$, and for any measurable $g:\left(0,r\right)\rightarrow[0,\infty]$  we have the inequality:
$$\int\limits_0^r g\left(\rho\right)\left(\int\limits_{S_\rho\left(a\right)}\left\lvert df_\rho\right\rvert^nd\sigma_{S_\rho\left(a\right)}\right)d\rho\leq \int\limits_{B_r\left(a\right)}g\left(\left\lvert x-a\right\rvert\right)\left\lvert\nabla f\left(x\right)\right\rvert^ndx$$
Where $\sigma_{S_\rho\left(a\right)}$ is the surface measure, and $df_\rho$ is the differential of the trace function $\trace_{S_\rho\left(a\right)} f$. 
\end{claim}
\begin{proof}
By Claim \ref{rep tr} we can assume $f=\hat{f}_a$. Take some sequence of smooth functions $\left\{f^k\right\}_{k\in\mathbb{N}}$ converging to $f$ in $W^{1,n}\left(B_r\left(a\right)\right)$. Then if we consider spherical coordinates $\left(\rho,\theta\right)$, and define $\frac{\partial f}{\partial\theta}$ as the matrix obtained from $\nabla f$ by removing the first column (we know that same relation holds for strong derivatives, thus all $f^k$), by the Fubini-Tonelli Theorem, it holds that:
\begin{equation*}
    \begin{split}
        & \int\limits_0^r\left(\int\limits_{S_\rho\left(a\right)}\left\lvert f^k-f\right\rvert^n+\left\lvert \frac{\partial f^k}{\partial\theta}-\frac{\partial f}{\partial\theta}\right\rvert^nd\sigma_{S_\rho\left(a\right)}\right)d\rho\\
        &\qquad =\int\limits_{B_r\left(a\right)}\left\lvert f^k-f\right\rvert^n+\left\lvert \frac{\partial f^k}{\partial\theta}-\frac{\partial f}{\partial\theta}\right\rvert^ndx\leq\int\limits_{B_r\left(a\right)}\left\lvert f^k-f\right\rvert^n+\left\lvert\nabla f^k-\nabla f\right\rvert^ndx
    \end{split}
\end{equation*}
Which tends to zero by the assumed $W^{1,n}\left(B_r\left(a\right)\right)$ convergence. Hence, for a.e.\ $\rho\in(0,r_a)$ we have:
$$\lim\limits_{k\rightarrow\infty}\int\limits_{S_\rho\left(a\right)}\left\lvert f^k-f\right\rvert^n+\left\lvert \frac{\partial f^k}{\partial\theta}-\frac{\partial f}{\partial\theta}\right\rvert^nd\sigma_{S_\rho\left(a\right)}=0$$
Which means that $f\in W^{1,n}\left(S_\rho\left(a\right)\right)$, with $df_\rho=\frac{\partial f}{\partial\theta}$. The inequality follows from applying the Fubini-Tonelli Theorem and noting that by construction:
$$\left\lvert df_\rho\right\rvert=\left\lvert\frac{\partial f}{\partial\theta}\right\rvert\leq\left\lvert\nabla f\right\rvert.$$
\end{proof}
\begin{definition}
\label{defF}
Let $a\in\Omega$, then we define the set:
$$F\left(a\right)=F\left(a,f\right)=\bigcap\limits_{\rho\in\left(0,r_a\right)\backslash Z_a} E\left(f,B_\rho\left(a\right)\right)$$
Where $Z_a=\left\{\rho\in\left(0,r_a\right)\mid f\notin W^{1,n}\left(S_\rho\left(a\right)\right)\right\}\cup Z_{a,f}$. \\
Additionally if $A\subset\Omega$, define $F\left(A\right)=\bigcup\limits_{a\in A}F\left(a\right)$.
\end{definition}
\begin{remark}
By Lemma \ref{mono E} $F\left(a\right)$ is well defined, non-empty, and compact as the the intersection of a descending chain of non-empty compact sets.
\end{remark}
\subsection{Proof of Regularity by Degree}
The first step to prove the regularity result, is to observe the following Lemma, concerning the last definition made in the previous section. (First two properties are important to prove the result, the rest are for completeness).
\begin{lemma} [properties of $F\left(a\right)$]
\label{prop of F}
Let $a\in\Omega$, then we have the following:
\begin{enumerate}
    \item For every $\rho\in\left(0,r_a\right)\backslash Z_a$ we have that $\diam E\left(f,B_\rho\left(a\right)\right)=\underset{S_\rho\left(a\right)}{\osc}\left(f\right)$, and in particular: 
    $$\diam F\left(a\right)=\lim\limits_{\underset{\rho\in\left(0,r_a\right)\backslash Z_a}{\rho\rightarrow0}}{\diam E\left(f,B_\rho\left(a\right)\right)}=\lim\limits_{\underset{\rho\in\left(0,r_a\right)\backslash Z_a}{\rho\rightarrow0}}{\underset{S_\rho\left(a\right)}{\osc}f}.$$
    \item Let $\tilde{f}=\left(\tilde{f}^1,...,\tilde{f}^n\right)$ be the representative of $f=\left(f^1,...,f^n\right)$ defined by:
    \begin{equation}
    \label{lebesgue rep}
        \tilde{f}^i\left(x\right)=\limsup\limits_{\rho\rightarrow0^+}{\frac{1}{\lambda_n\left(B_\rho\left(x\right)\right)}\int\limits_{B_\rho\left(x\right)}f^i\left(z\right)dz}.
    \end{equation}
    Then we have the inequality:
    $$\limsup\limits_{\rho\rightarrow0^+}{\underset{B_\rho\left(x\right)}{\osc}\tilde{f}}\leq\diam F\left(x\right),$$
    and in particular $\tilde{f}$ is continuous at any $x\in\Omega$ such that $diam F\left(x\right)=0$.
    \item $F\left(a\right)$ is connected.
    \item If $f$ has a representative $\bar{f}$ continuous at $a$ then $\diam F\left(a\right)=0$, and in particular $F\left(a\right)=\left\{\bar{f}\left(a\right)\right\}$.
    \item Let $\Omega_1\subset\overline{\Omega_1}\subset\Omega$ open with smooth boundary, and denote $\Gamma_1=\partial\Omega_1$.\\
    Then if
    $f\in W^{1,n}\left(\Gamma_1\right)$, with continuous representative $\bar{f}$ of $\trace_{\Gamma_1} f$, then $F\left(x\right)\subset E\left(f,\Omega_1\right)$, and $\bar{f}\left(x\right)\in F\left(x\right)$ for every $x\in\Gamma_1$.
    \item For every $V\subset\mathbb{R}^n$ open, the set $\left\{x\in \Omega\mid F\left(x\right)\subset V\right\}$ is also open. Moreover the function $x\rightarrow \diam F\left(x\right)$ is upper-semicontinuous. 
    \item For all $K\subset\Omega$ compact, $F\left(K\right)$ is also compact.
    \item Let $\Omega_1$ as in point five, then $E\left(f,\Omega_1\right)\subset F\left(\bar{\Omega}_1\right)$.
\end{enumerate}
\end{lemma}
\begin{proof}
\begin{enumerate}
    \item Note that as $E\left(f,B_\rho\left(a\right)\right)$ is the union of connected components (without the outermost component) of $\mathbb{R}^n\backslash f\left(S_\rho\left(a\right)\right)$ and $f\left(S_\rho\left(a\right)\right)$, we have:
    $$\partial E\left(f,B_\rho\left(a\right)\right)\subset f\left(S_\rho\left(a\right)\right)\subset E\left(f,B_\rho\left(a\right)\right)$$ 
    And using the compactness of $E\left(f,B_\rho\left(a\right)\right)$, we can observe that:
    \begin{equation*}
        \begin{split}
            \diam E\left(f,B_\rho\left(a\right)\right) & =\diam\partial E\left(f,B_\rho\left(a\right)\right)\leq\diam f\left(S_\rho\left(a\right)\right) \\
            & \leq\diam E\left(f,B_\rho\left(a\right)\right)
        \end{split}
    \end{equation*}
    Thus we have that 
    $$\diam E\left(f,B_\rho\left(a\right)\right)=\diam f\left(S_\rho\left(a\right)\right)=\underset{S_\rho\left(a\right)}{\osc}\left(f\right)$$
    the second part follows directly from the definition of $F$ and Lemma \ref{mono E}.
    \item Let $x\in\Omega$, and consider $\rho\in\left(0,r_x\right)\backslash Z_x$ we will show that 
    $$\underset{B_\rho\left(x\right)}{\osc}\tilde{f}\leq\diam E\left(f,B_\rho\left(x\right)\right)$$ 
    and the claim follows directly from that. The first step is to show that for every $ z_0\in B_\rho\left(x\right)$, and small enough $r>0$ we have that $f\left(z\right)\in E\left(f,B_\rho\left(x\right)\right)$ for a.e.\ $z\in B_r\left(z_0\right)$. Indeed, as $B_\rho\left(x\right)$ is open, we can choose a small enough $r_0>0$ such that $f\in W^{1,n}\left(S_{r_0}\left(z_0\right)\right)$ and also $B_{r_0}\left(z_0\right)\subset B_\rho\left(x\right)$, thus by Corollary \ref{belong} and Lemma \ref{mono E} we have that $f\left(z\right)\in E\left(f,B_\rho\left(x\right)\right)$ for a.e.\ $z\in B_{r_0}\left(z_0\right)$, but then the same is true for a.e.\ $0<r<r_0$.\\
    Now we can note that for $r>0$ small enough, each
    $$w_r\left(z_0\right)=\frac{1}{\lambda_n\left(B_r\left(z_0\right)\right)}\int\limits_{B_r\left(z_0\right)}f\left(z\right)dz$$ 
    is the average of elements inside $E\left(f,B_\rho\left(x\right)\right)$, thus we must have that $w_r\left(z_0\right)\in\convex E\left(f,B_\rho\left(x\right)\right)$, where $\convex A$ denotes the convex hull of the set $A$.  But as the convex hull of a compact set, $\convex E\left(f,B_\rho\left(x\right)\right)$ is compact as well, thus $\tilde{f}\left(z_0\right)\in\convex E\left(f,B_\rho\left(x\right)\right)$ as a limiting point of $w_r\left(z_0\right)$. Thus using the fact that taking a convex hull conserves the diameter of a set, we obtain that: 
    $$\underset{B_\rho\left(x\right)}{\osc}\tilde{f}=\diam\tilde{f}\left(B_\rho\left(x\right)\right)\leq\diam\convex E\left(f,B_\rho\left(x\right)\right)=\diam E\left(f,B_\rho\left(x\right)\right)$$
    As we desired.
    \item Assume by contradiction $F\left(a\right)$ is not connected, then there is a continuous function $g:F\left(a\right)\rightarrow\left\{0,1\right\}$, with $g\left(x_0\right)=0,g\left(x_1\right)=1$ for some $x_1,x_0\in F\left(a\right)$. By definition take any $\rho_1\in\left(0,\min\left\{r_a,1\right\}\right)\backslash Z_a$, and we have $F\left(a\right)\subset E\left(f,B_{\rho_1}\left(a\right)\right)$, and because $F\left(a\right)$ is compact, we can extend $g$ to some continuous $G:E\left(f,B_{\rho_1}\left(a\right)\right)\rightarrow[0,1]$, and by the Intermediate Value Theorem (uses that $E\left(f,B_{\rho_1}\left(a\right)\right)$ is connected) we have some $x_{\rho_1}\in E\left(f,B_{\rho_1}\left(a\right)\right)$ with $G\left(x_{\rho_1}\right)=\frac{1}{2}$. Now take inductively $\rho_n\in\left(0,\min\left\{\rho_{n-1},\frac{1}{n}\right\}\right)\backslash Z_a$ and by Lemma \ref{mono E}, we have $F\left(a\right)\subset E\left(f,B_{\rho_n}\left(a\right)\right)\subset E\left(f,B_{\rho_1}\left(a\right)\right)$ thus we can look at $\left.G\right\rvert_{E\left(f,B_{\rho_n}\left(a\right)\right)}$ and again find some $x_{\rho_n}\in E\left(f,B_{\rho_n}\left(a\right)\right)$ with $G\left(x_{\rho_n}\right)=\frac{1}{2}$. and by taking limits, we must have $x_{\frac{1}{2}}=\lim\limits_{n\rightarrow\infty}x_{\rho_n}\in F\left(a\right)$ but then
    $$g\left(x_{\frac{1}{2}}\right)=G\left(x_{\frac{1}{2}}\right)=G\left(\lim\limits_{n\rightarrow\infty}x_{\rho_n}\right)=\lim\limits_{n\rightarrow\infty}G\left(x_{\rho_n}\right)=\lim\limits_{n\rightarrow\infty}\frac{1}{2}=\frac{1}{2}$$
    Giving a contradiction, thus $F\left(a\right)$ is connected.
    \begin{remark}
    The argument here is a general one, which shows that the intersection of descending  chain of compact connected sets is also connected.
    \end{remark}
    \item By continuity, for every $n\in\mathbb{N}$ there is a $0<\delta_n<\frac{1}{n}$ satisfying 
    $$\bar{f}\left(B_{\delta_n}\left(a\right)\right)\subset B_\frac{1}{n}\left(\bar{f}\left(a\right)\right)$$
    By Corollary \ref{cor rep} we can take some $r_n\in\left(0,\min\left\{\frac{1}{n},\delta_n,r_a\right\}\right)\backslash \left(Z_a\cup Z_{a,\bar{f}}\right)$ (Recall that $Z_a$ is defined with respect to $f$ and not $\Bar{f}$). Note that for almost every $x\in S_{r_n}\left(a\right)$ we have $f\left(x\right)=\bar{f}\left(x\right)\in B_\frac{1}{n}\left(\bar{f}\left(a\right)\right)$. Additionally by continuity of $f$ on $S_{r_n}\left(a\right)$ we get the inclusion $f\left(S_{r_n}\left(a\right)\right)\subset\overline{B_\frac{1}{n}\left(\bar{f}\left(a\right)\right)}$. Using this and the first point of this Lemma we obtain:
    $$\diam F\left(a\right)=\lim\limits_{\underset{\rho\in\left(0,r_a\right)\backslash Z_a}{\rho\rightarrow0}}{\underset{S_\rho\left(a\right)}{\osc}f}=\lim\limits_{n\rightarrow\infty}{\underset{S_{r_n}\left(a\right)}{\osc}f}\leq\lim\limits_{n\rightarrow\infty}\frac{2}{n}=0$$
    Thus $\diam F\left(a\right)=0$. Because of this, $F\left(a\right)$ is a singleton, hence, to prove that $F\left(a\right)=\left\{\bar{f}\left(a\right)\right\}$ it suffices to just show $\bar{f}\left(a\right)\in F\left(a\right)$.\\
    As before, take $r_n$ and consider any $z_n\in S_{r_n}\left(a\right)$ satisfying the inclusion $\bar{f}\left(z_n\right)=f\left(z_n\right)\in E\left(f,B_{r_n}\left(a\right)\right)$ (the belong relation is given by the definition of $E\left(f,B_{r_n}\left(a\right)\right)$). Then we can observe that $\bar{f}\left(z_n\right)\in B_\frac{1}{n}\left(\bar{f}\left(a\right)\right)\cap E\left(f,B_{r_n}\left(a\right)\right)$ and by taking limits, we obtain $\bar{f}\left(a\right)=\lim\limits_{n\rightarrow\infty}\bar{f}\left(z_n\right)\in F\left(a\right)$ (as $r_n\underset{n\rightarrow\infty}{\rightarrow}0$) as needed.
    \item Because $\Omega_1$ is open, for any $x\in\Omega_1$ consider some $r\in \left(0,r_x\right)\backslash Z_x$ with $B_r\left(x\right)\subset\Omega_1$ then by Lemma \ref{mono E} we have
    $$F\left(x\right)=\bigcap\limits_{\rho\in\left(0,r_x\right)\backslash Z_x} E\left(f,B_\rho\left(x\right)\right)\subset E\left(f,B_r\left(x\right)\right)\subset E\left(f,\Omega_1\right)$$
    Now let $x\in\Gamma_1$, and $\left\{r_n\right\}_{n\in\mathbb{N}}\subset\left(0,r_x\right)\backslash Z_x$ be a sequence decreasing to zero, and consider for each $n$ some small segment, around $x$, of $\Gamma_1$, which is contained on $B_{r_n}\left(x\right)$, then as in the proof of Lemma \ref{mono E}, we have $\bar{f}\left(x\right)\in E\left(f,B_{r_n}\left(x\right)\right)$, and this is true for every $n\in\mathbb{N}$. Hence, $\bar{f}\left(x\right)\in F\left(x\right)$ as needed.
    \item Let $x\in\Omega$ s.t.\ $F\left(x\right)\subset V$. Let's first see that there is some $r\in\left(0,r_x\right)\backslash Z_x$ with $E\left(f,B_r\left(x\right)\right)\subset V$.
    Assume the contrary, and for every $n\in\mathbb{N}$ take some $r_n\in \left(0,\min\left\{r_{n-1},\frac{1}{n}\right\}\right)\backslash Z_x$ (take $r_0=r_x$), and consider $y_n\in E\left(f,B_{r_n}\left(x\right)\right)\backslash V$. But then using the definition of $F\left(x\right)$, Lemma \ref{mono E}, the openness of $V$, and compactness of $E\left(f,B_{r_1}\left(x\right)\right)$, we can pass to an appropriate subsequence to obtain:
    $$\lim\limits_{n\rightarrow\infty}y_n\in F\left(x\right)\backslash V=\emptyset$$
    Giving a contradiction.\\
    Then for every $z\in B_\frac{r}{3}\left(x\right)$ take some $r'\in\left(0,\min\left\{r_z,\frac{r}{3}\right\}\right)\backslash Z_z$ and again by definition of $F\left(z\right)$ and Lemma \ref{mono E} we have:
    $$F\left(z\right)=\bigcap\limits_{\rho\in\left(0,r_z\right)\backslash Z_z} E\left(f,B_\rho\left(z\right)\right)\subset E\left(f,B_{r'}\left(z\right)\right)\subset E\left(f,B_r\left(x\right)\right)\subset V$$
    Thus the set $\left\{x\in \Omega\mid F\left(x\right)\subset V\right\}$ is open.\\
    The upper-semicontinuity Claim follows from  proving that for every $\epsilon>0$ the set:
    $$\left\{x\in\Omega\mid\diam F\left(x\right)<\epsilon\right\}$$
    is open.
    Indeed if $x\in\Omega$ is such that $d=\diam F\left(x\right)<\epsilon$, then  consider $\tilde{\epsilon}=\frac{\epsilon-d}{2}>0$, and look at the set:
    $$\tilde{F}\left(x\right)=\left\{y\in\Omega\mid\dist\left(y,F\left(x\right)\right)<\tilde{\epsilon}\right\}$$
    It easily follows that this set is open and satisfies that $F\left(x\right)\subset\tilde{F}\left(x\right)$, thus, by the previous part, the set:
    $$\left\{y\in\Omega\mid F\left(y\right)\subset\tilde{F}\left(x\right)\right\}$$
    Is also open, and has that, by Triangle Inequality, for every $y$ in it the bound:
    $$\diam F\left(y\right)\leq\diam\tilde{F}\left(x\right)\leq\diam F\left(x\right)+2\tilde{\epsilon}<d+2\cdot\frac{\epsilon-d}{2}=\epsilon$$
    Is satisfied, giving the desired result.
    \item Let $\left\{y_n\right\}_{n\in\mathbb{N}}\in F\left(K\right)$. Then for each $n\in\mathbb{N}$ there is some $x_n\in K$ with $y_n\in F\left(x_n\right)$, by compactness of $K$, we can take a subsequence of $\left\{x_n\right\}_{n\in\mathbb{N}}$ converging to some $x\in K$. Take some further subsequence such that $\dist\left(x,x_n\right)$ is monotone decreasing.  And take some decreasing sequence of radii $\left\{r_n\right\}_{n\in\mathbb{N}}\subset \left(0,r_x\right)\backslash Z_x$ with $r_n\underset{n\rightarrow\infty}{\rightarrow}0$ and $x_n\in B_{r_n}\left(x\right)$ (assume subsequence starts at large enough $n$ such that the second requirement can be achieved). Now note that for every $n\in\mathbb{N}$ we have some $r'_n\in \left(0,r_{x_n}\right)\backslash Z_{x_n}$ satisfying $B_{r'_n}\left(x_n\right)\subset B_{r_n}\left(x\right)$. Thus by definitions and Lemma \ref{mono E} we finally have: 
    \begin{equation*}
        \begin{split}
            y_n\in F\left(x_n\right)=\bigcap\limits_{\rho\in\left(0,r_{x_n}\right)\backslash Z_{x_n}} E\left(f,B_\rho\left(x_n\right)\right) & \subset E\left(f,B_{r'_n}\left(x_n\right)\right)\\
            & \subset E\left(f,B_{r_n}\left(x\right)\right)
        \end{split}
    \end{equation*}
    Thus using Lemma \ref{mono E} and compactness of $E\left(f,B_{r_1}\left(x\right)\right)$, we can take a further subsequence of $\left\{y_n\right\}_{n\in\mathbb{N}}$ which converges. Since $\left\{y_n\right\}_{n\geq N}\subset E\left(f,B_{r_N}\left(x\right)\right)$ we have $\lim\limits_{n\rightarrow\infty}y_n\in B_{r_n}\left(x\right)$ for every $N\in\mathbb{N}$ hence $\lim\limits_{n\rightarrow\infty}y_n\in F\left(x\right)\subset F\left(K\right)$. Thus $F\left(K\right)$ is a compact set.
    \item Note that from the fifth point we get $\bar{f}\left(\Gamma_1\right)\subset F\left(\Gamma_1\right)\subset F\left(\overline{\Omega_1}\right)$.\\ 
    Thus we observe that the set
    \begin{equation*}
        \begin{split}
            H   & =E\left(f,\Omega_1\right)\backslash F\left(\overline{\Omega_1}\right)=E\left(f,\Omega_1\right)\backslash\left(\bar{f}\left(\Gamma_1\right)\cup F\left(\overline{\Omega_1}\right)\right)
    \\& =\left(E\left(f,\Omega_1\right)\backslash\bar{f}\left(\Gamma_1\right)\right)\backslash F\left(\overline{\Omega_1}\right)
        \end{split}
    \end{equation*}
    is open because (i) $E\left(f,\Omega_1\right)\backslash\bar{f}\left(\Gamma_1\right)$ is, by definition, the union of connected components of $\mathbb{R}^n\backslash\bar{f}\left(\Gamma_1\right)$, and thus open. And, (ii) by the seventh point, we know that the set $F\left(\overline{\Omega_1}\right)$ is compact, and in particular, closed.\\
    Assume by contradiction that $H\neq\emptyset$ and let $g\in C^\infty$ with support inside a closed ball $B\subset H$ and also satisfying $\int\limits_{\mathbb{R}^n} g\left(z\right)dz=1$.\\
    Let $x\in\overline{\Omega_1}$. Because $F\left(x\right)\cap B=\emptyset$, then using the claim made in proof of the sixth part (use that $B^C$ is open), there is some $r\in \left(0,r_x\right)\backslash Z_x$ with $E\left(f,B_r\left(x\right)\right)\cap B=\emptyset$. Thus by Corollary \ref{belong} we know that $f\left(z\right)\in\mathbb{R}^n\backslash B$ of almost every $z\in B_r\left(x\right)$. Using the compactness of $\overline{\Omega_1}$ we obtain the same result for almost every $z\in\overline{\Omega_1}$. Thus, by Formula \eqref{change of var}, for every $y\in B$ we have:
    $$\deg \left(f,\Omega_1,y\right)=\deg \left(f,\Omega_1,y\right)\int\limits_{\mathbb{R}^n} g\left(z\right)dz=\int\limits_{\Omega_1} g\left(f\left(x\right)\right) \det \nabla f\left(x\right)dx=0$$
    On the other hand, because $y\in E\left(f,\Omega_1\right)$, by definition we must have $\deg \left(f,\Omega_1,y\right)\geq1$, giving a contradiction. Thus $H=\emptyset$.
\end{enumerate}
\end{proof}
After the following Lemma we are finally ready to prove the desired regularity result.
\begin{lemma}[Oscillation Lemma]\cite[Lemma 4.10.1]{iwaniec}
\label{osc lem}
Let $a\in\Omega$ and $r\in\left(0,r_a\right)\backslash Z_a$, then the following inequality holds:
$$\left(\underset{S_r\left(a\right)}{\osc}f\right)^n\leq C\cdot r\int\limits_{S_r\left(a\right)}\left\lvert df\right\rvert^nd\sigma _{S_r\left(a\right)}$$
Where $\sigma_{S_\rho\left(a\right)}$ is the surface measure, $f$ is the appropriate continuous representative of the trace function $\trace_{S_\rho\left(a\right)} f$, $df$ is the differential of that representative and $C$ is some positive constant depending on $n$. 
\end{lemma}
\begin{remark}
This is essentially Morrey's Inequality for functions on the sphere.
\end{remark}
\begin{remark}
On \cite{iwaniec} the Lemma is stated in a more general way, but for Lipshitz functions, but by a density argument, it can be extended to our functions, since $W^{1,n}\left(S_r\left(a\right)\right)$ convergence implies uniform convergence (uses $\dim S_r\left(a\right)<n$).
\end{remark}
\begin{restatable}{thm}{regresult}
Let $f\in W^{1,n}\left(\Omega\right)$ such that $\det\nabla f>0$ a.e.\ in $\Omega$ then the representative $\tilde{f}$ of $f$, given by Formula \eqref{lebesgue rep}, is continuous.
\end{restatable}
\begin{proof}
Let's start by defining the set
$$S=\left\{x\in\Omega\mid\diam F\left(x\right)>0\right\}$$
By the second point of Lemma \ref{prop of F} we get that $\tilde{f}$, the representative given by the Lebesgue Differentiation Theorem, is continuous on the complement of $S$. To prove the Theorem we just need to show that $S$ is empty.\\
Assume that $x\in S$, then by the first part of Lemma \ref{prop of F}, and Lemma \ref{osc lem} we have:
\begin{equation*}
    \begin{split}
        0 & <\left(\diam F\left(x\right)\right)^n=\left(\lim\limits_{\underset{\rho\in\left(0,r_x\right)\backslash Z_x}{\rho\rightarrow0}}{\underset{S_\rho\left(x\right)}{\osc}f}\right)^n\\
 & =\lim\limits_{\underset{\rho\in\left(0,r_x\right)\backslash Z_x}{\rho\rightarrow0}}\left({\underset{S_\rho\left(x\right)}{\osc}f}\right)^n<\lim\limits_{\underset{\rho\in\left(0,r_x\right)\backslash Z_x}{\rho\rightarrow0}}C\cdot \rho\int\limits_{S_\rho\left(x\right)}\left\lvert df\right\rvert^nd\sigma _{S_\rho\left(x\right)}
    \end{split}
\end{equation*}
The last inequality implies the following growth inequality, for $\rho\in \left(0,r_x\right)\backslash Z_x$ small enough
$$\int\limits_{S_\rho\left(x\right)}\left\lvert df\right\rvert^nd\sigma _{S_\rho\left(x\right)}>\frac{C}{\rho}$$
For some new constant $C$. Thus we have the inclusion
$$S\subset\left\{x\in\Omega\mid\int\limits_0^{r_x}\left(\int\limits_{S_\rho\left(x\right)}\left\lvert df\right\rvert^nd\sigma _{S_\rho\left(x\right)}\right)d\rho=\infty\right\}$$
Now use Claim \ref{ineq of surf} with the constant function $g\left(x\right)=1$ we obtain
$$\int\limits_0^{r_x}\left(\int\limits_{S_\rho\left(x\right)}\left\lvert df\right\rvert^nd\sigma _{S_\rho\left(x\right)}\right)d\rho\leq\int\limits_{B_{r_x}\left(x\right)}\left\lvert\nabla f\left(z\right)\right\rvert^ndz$$
This gives the further inclusion 
$$S\subset\left\{x\in\Omega\mid\int\limits_{B_{r_x}\left(x\right)}\left\lvert\nabla f\left(z\right)\right\rvert^ndz=\infty\right\}$$
But the right hand set is empty as $f\in W^{1,n}\left(B_{r_x}\left(x\right)\right)$ for any $x\in\Omega$. Thus the set $S$ is empty meaning that the representative $\tilde{f}$ is continuous in all of $\Omega$.
\end{proof}
\section{Proof by Oscillation Argument}
\subsection{Idea}
This method is based on proving an inequality bounding the oscillation of the function on a ball by the oscillation on just the surface of that ball, and subsequently, exploiting the already known inequalities of Sobolev functions to obtain the desired regularity result. In doing so a further general statement of regularity of integrable functions must be made.
\subsection{Oscillation Inequality}
The first notion needed for this proof is that of essential oscillation:
\begin{definition}[Essential Oscillation]
\label{eosc}
Given a measure space $\left(X,\mathcal{M},\mu\right)$, a metric space $\left(Y,d\right)$, and a function $f:X\rightarrow Y$ we define the essential oscillation of $f$ on $B\subset X$ by:
$$\underset{B}{\eosc}f=\inf\left\{\diam f(A)\mid A\subset B,\text{ }\mu\left(B\backslash A\right)=0\right\}$$
\end{definition}
\begin{remark}
Whenever $f$ is continuous, and $\mu$ is Radon, the essential oscillation agrees with the classical one.
\end{remark}
Using this definition we can make the following Claims: 
\begin{claim}
\label{ball}
Let $f\in W^{1,n}\left(\Omega;\mathbb{R}^n\right)$ with $\det\nabla f\geq0$ a.e.\ Let $x\in\Omega$, and $r\in\left(0,r_x\right)\backslash Z_x$, (for definition of $r_x$ and $Z_x$, see Chapter 2, Definition \ref{defF}). If we take a ball $\overline{B}$ containing $f\left(S_r\left(x\right)\right)$, we have that: 
$$\lambda_n\left(f^{-1}\left(\mathbb{R}^n\backslash\overline{B}\right)\cap B_r\left(x\right)\right)=0,$$
which can be interpreted as $f\left(y\right)\in\overline{B}$ for a.e.\ $y\in B_r\left(x\right)$.
\end{claim}
\begin{remark}
    The assumption $r\notin Z_a$ is to ensure the existence of a ball $\overline{B}$ with $f\left(S_r\left(x\right)\right)\subset\overline{B}$.
\end{remark}
\begin{proof}
Assume by contradiction that $K=f^{-1}\left(\mathbb{R}^n\backslash \overline{B}\right)\cap B_r\left(x\right)$ has positive measure. First, note that as $J_f=\det\nabla f>0$ a.e., we have that:
$$\int\limits_KJ_fdy>0$$
Now consider the radial retract $\rho:\mathbb{R}^n\rightarrow\overline{B}$, and since $f\left(K\right)\subset\mathbb{R}^n\backslash\overline{B}$, we have that $\rho\circ f$, sends $K$, a set with positive $n$-dimensional Lebesgue measure, to $\partial\Bar{B}$, an $\left(n-1\right)$-dimensional set, which implies that $\left.J_{\rho\circ f}\right\rvert_K=0$, hence:
$$\int\limits_KJ_fdy>0=\int\limits_KJ_{\rho\circ f}dy$$
Also recall that by definition $\left.\rho\right\rvert_{\overline{B}}=\identity_{\overline{B}}$, and note that if we denote $H=B_r\left(x\right)\backslash K$, then $f\left(H\right)\subset\overline{B}$, thus $\left.\rho\circ f\right\rvert_H=\left.\rho\right\rvert_{f\left(H\right)}\circ\left.f\right\rvert_H=\left.f\right\rvert_H$, hence:
$$\int\limits_HJ_fdy=\int\limits_HJ_{\rho\circ f}dy$$
On the other hand, we recall that we defined $\overline{B}$ such that $f\left(S_r\left(x\right)\right)\subset\overline{B}$, and thus we also have that $\left.\rho\circ f\right\rvert_{S_r\left(x\right)}=\left.\rho\right\rvert_{f\left(S_r\left(x\right)\right)}\circ\left.f\right\rvert_{S_r\left(x\right)}=\left.f\right\rvert_{S_r\left(x\right)}$, thus using the fact that the Jacobian is a Null-Lagrangian, we obtain that:
\begin{equation*}
\begin{split}
    \int\limits_{B_r\left(x\right)}J_{\rho\circ f}dy=\int\limits_{B_r\left(x\right)}J_fdy=\int\limits_KJ_fdy+\int\limits_HJ_fdy & >\int\limits_KJ_{\rho\circ f}dy+\int\limits_HJ_{\rho\circ f}dy\\ & =\int\limits_{B_r\left(x\right)}J_{\rho\circ f}dy,
\end{split}
\end{equation*}
which is a clear contradiction.
\end{proof}
\begin{claim}
\label{osc ineq}
Given a map $f\in W^{1,n}\left(\Omega;\mathbb{R}^n\right)$ with $\det\nabla f\geq0$ a.e.\, then for every $x\in\Omega$, and every $r\in\left(0,r_x\right)\backslash Z_x$, the following inequality holds:
\begin{equation}
\underset{B_r\left(x\right)}{\eosc}f\leq2\underset{S_r\left(x\right)}{\osc}f . 
\end{equation}
\end{claim}
\begin{proof}
Consider a ball $\overline{B}$ with $f\left(S_r\left(x\right)\right)\subset\overline{B}$, and $\diam \overline{B}\leq 2\diam f\left(S_r\left(x\right)\right)$. Then, according to Claim \ref{ball}, if we consider $A=\left\{y\in\ B_r\left(x\right)\mid f\left(y\right)\in\overline{B}\right\}\subset B_r\left(x\right)$, it has $\lambda_n\left(B_r\left(x\right)\backslash A\right)=0$, thus:
$$\underset{B_r\left(x\right)}{\eosc}f\leq\diam f\left(A\right)\leq\diam\overline{B}\leq2\diam f\left(S_r\left(x\right)\right)=2\underset{S_r\left(x\right)}{\osc}f$$
as needed.
\end{proof}
\subsection{Proof of Regularity}
The first step of this proof is to identify the canonical continuous representative of integrable functions.
\begin{theorem}[Canonical Continuous Representative]
\label{canrep}
Let $f\in L^1\left(\Omega\right)$. Then the representative $\tilde{f}$ of $f$ given by the Lebesgue Differentiation Theorem is continuous at $x\in\Omega$ whenever the equality $\lim\limits_{r\rightarrow0}\underset{B_r\left(x\right)}{\eosc}f=0$ holds.
\end{theorem}
\begin{remark}
This theorem means that an $L^1$ function has a continuous representative at any point $x$ if and only if the representative given by the Lebesgue Differentiation Theorem is continuous at $x$.
\end{remark}
\begin{proof}
Let $x\in\Omega$ such that $\lim\limits_{r\rightarrow0}\underset{B_r\left(x\right)}{\eosc}f=0$, and assume by contradiction that $\tilde{f}$ is not continuous at $x$. Then we have some $\epsilon>0$ such that for any $\delta>0$ there is some $x'_\delta\in B_\delta\left(x\right)$ (from now on just $x'$) satisfying $\left\lvert\tilde{f}\left(x\right)-\tilde{f}\left(x'\right)\right\rvert\geq\epsilon$.\\
Using the definition of limit, we have some $\rho>0$ with $\underset{B_\rho\left(x\right)}{\eosc}f<\frac{\epsilon}{5}$, and as above, there is some $x'\in B_\rho\left(x\right)$ with $\left\lvert\tilde{f}\left(x\right)-\tilde{f}\left(x'\right)\right\rvert\geq\epsilon$.\\
By definition of essential oscillation, we can consider a set $A\subset B_\rho\left(x\right)$ of full measure satisfying $\diam f\left(A\right)<\frac{\epsilon}{5}$.\\
By definition of $\tilde{f}\left(x\right)$ there is some $\rho>r_1>0$ such that:
$$\left\lvert\frac{1}{\lambda_n\left(B_{r_1}\left(x\right)\right)}\int\limits_{B_{r_1}\left(x\right)}f\left(z\right)dz-\tilde{f}\left(x\right)\right\rvert<\frac{\epsilon}{5}.$$
We notice that: 
\begin{equation*}
    \begin{split}
        &\frac{1}{\lambda_n\left(B_{r_1}\left(x\right)\right)}\int\limits_{B_{r_1}\left(x\right)}f\left(z\right)dz\\
        &\qquad =\frac{1}{\lambda_n\left(B_{r_1}\left(x\right)\cap A\right)}\int\limits_{B_{r_1}\left(x\right)\cap A}f\left(z\right)dz\in\convex f\left(B_{r_1}\left(x\right)\cap A\right).
    \end{split}
\end{equation*}
Since the diameter of the convex hull of a set equals to the diameter of the set,
$$\diam\convex f\left(B_{r_1}\left(x\right)\cap A\right)=\diam f\left(B_{r_1}\left(x\right)\cap A\right)\leq \diam f\left(A\right)<\frac{\epsilon}{5}.$$
Thus, by the definitions of diameter and convex hull, and the triangle inequality, for every $y\in B_{r_1}\left(x\right)\cap A$ we have that:
\begin{equation*}
    \begin{split}
        \left\lvert f\left(y\right)-\tilde{f}\left(x\right)\right\rvert
& \leq\left\lvert f\left(y\right)-\frac{1}{\lambda_n\left(B_{r_1}\left(x\right)\right)}\int\limits_{B_{r_1}\left(x\right)}f\left(z\right)dz\right\rvert\\
& + \left\lvert\frac{1}{\lambda_n\left(B_{r_1}\left(x\right)\right)}\int\limits_{B_{r_1}\left(x\right)}f\left(z\right)dz-\tilde{f}\left(x\right)\right\rvert\\
& <\frac{\epsilon}{5}+\frac{\epsilon}{5}=\frac{2\epsilon}{5}
    \end{split}
\end{equation*}
In a similar manner, we choose some $r_2>0$ such that $B_{r_2}\left(x'\right)\subset B_\rho\left(x\right)$ and also:
$$\left\lvert\frac{1}{\lambda_n\left(B_{r_2}\left(x'\right)\right)}\int\limits_{B_{r_2}\left(x'\right)}f\left(z\right)dz-\tilde{f}\left(x'\right)\right\rvert<\frac{\epsilon}{5}$$
And following the same process as above, we see that for any $y\in B_{r_2}\left(x'\right)\cap A$, the inequality $\left\lvert f\left(y\right)-\tilde{f}\left(x'\right)\right\rvert<\frac{2\epsilon}{5}$ also holds.\\
Now choose some $y_1\in B_{r_1}\left(x\right)\cap A$, and by the Inverse Triangle Inequality we see that:
$$\left\lvert f\left(y_1\right)-\tilde{f}\left(x'\right)\right\rvert\geq\left\lvert\left\lvert f\left(y_1\right)-\tilde{f}\left(x\right)\right\rvert-\left\lvert \tilde{f}\left(x'\right)-\tilde{f}\left(x\right)\right\rvert\right\rvert>\epsilon-\frac{2\epsilon}{5}=\frac{3\epsilon}{5}$$
If we also choose some $y_2\in B_{r_2}\left(x'\right)\cap A$, we note that $y_1,y_2\in A$, and by the Inverse Triangle Inequality we see that:
\begin{equation*}
    \begin{split}
        \diam f\left(A\right) & \geq\left\lvert f\left(y_1\right)-f\left(y_2\right)\right\rvert\geq\left\lvert\left\lvert f\left(y_1\right)-\tilde{f}\left(x'\right)\right\rvert-\left\lvert f\left(y_2\right)-\tilde{f}\left(x'\right)\right\rvert\right\rvert
        \\
        & >\frac{3\epsilon}{6}-\frac{2\epsilon}{5}=\frac{\epsilon}{5},
    \end{split}
\end{equation*}
which is a contradiction to the selection of $A$.
\end{proof}
Finally we proceed to prove the desired result:
\regresult*
\begin{proof}
By Theorem \ref{canrep}, it is only necessary to show that at any $x\in\Omega$ we have the limit $$\lim\limits_{r\rightarrow0}\underset{B_r\left(x\right)}{\eosc}f=0.$$
Indeed, let any $r\in\left(0,\frac{r_x}{2}\right)$, by Claims \ref{osc ineq}, and \ref{ineq of surf} (with $g\equiv1$), and Lemma \ref{osc lem} we see that: 
\begin{equation*}
    \begin{split}
        \left(\underset{B_r\left(x\right)}{\eosc}f\right)^n & =\frac{1}{\log2}\int\limits_r^{2r}\frac{\left(\underset{B_r\left(x\right)}{\eosc}f\right)^n}{\rho}d\rho\leq\frac{1}{\log2}\int\limits_r^{2r}\frac{\left(\underset{B_\rho\left(x\right)}{\eosc}f\right)^n}{\rho}d\rho\\
        & \overset{\ref{osc ineq}}{\leq}\frac{2^n}{\log2}\int\limits_r^{2r}\frac{\left(\underset{S_\rho\left(x\right)}{\osc}f\right)^n}{\rho}d\rho\overset{\ref{osc lem}}{\leq}\frac{2^n\cdot C}{\log2}\int\limits_r^{2r}\left(\int\limits_{S_\rho\left(x\right)}\left\lvert df\right\rvert^nd\sigma_{S_\rho\left(x\right)}\right)d\rho\\
        & \leq\frac{2^n\cdot C}{\log2}\int\limits_0^{2r}\left(\int\limits_{S_\rho\left(x\right)}\left\lvert df\right\rvert^nd\sigma_{S_\rho\left(x\right)}\right)d\rho\overset{\ref{ineq of surf}}{\leq}\frac{2^n\cdot C}{\log2}\int\limits_{B_{2r}\left(x\right)}\left\lvert \nabla f\right\rvert^ndz.
    \end{split}
\end{equation*}
Since $f\in W^{1,n}\left(\Omega\right)$, the last term tends to zero as $r\rightarrow0$, giving the desired limit.
\end{proof}
\section{Further Results - $VMO$ Theory}
\subsection{Objectives}
Until now we assumed that each function was at least $W^{1,n}$ but if we consider the steps of Chapter 2, we can ask if the same procedure applies to any larger family of functions that still have a defined degree. A larger class of functions that admits a definition  of degree (for a domain with boundary) is that of $VMO$ functions with $VMO$ trace. Thus in this Chapter we will construct the $VMO$ space, define the trace for the functions that admit it, and finally define the degree for such functions. Additionally, we will show that the function $\Tilde{F}\left(x,t\right)$ defined as in the introduction indeed belongs to such space.
\subsection{Definition of $VMO$ space}
The first step to understand $VMO$ (Vanishing Mean Oscillation) functions is to understand a larger space called $BMO$ (Bounded Mean Oscillation).
\begin{definition}[$BMO$ Space]
    Given any function $f\in L_\text{loc}^1\left(\Omega\right)$ we consider the semi-norm
    \begin{equation}
    \left\lVert f\right\rVert_{BMO} = \sup\limits_{\underset{x\in\Omega}{\epsilon<r_x}}\fint\limits_{B_\epsilon\left(x\right)}\left\lvert f\left(y\right)-\Bar{f}_\epsilon\left(x\right)\right\rvert dy
    \end{equation}
    where $r_x=\dist\left(x,\Gamma\right)$, $\fint\limits_A f\left(y\right)dy=\frac{1}{\lambda_n\left(A\right)}\int\limits_A f\left(y\right)dy$ ,and $\Bar{f}_\epsilon$ is given by
$$\Bar{f}_\epsilon\left(x\right)=\fint\limits_{B_\epsilon\left(x\right)}f\left(y\right)dy.$$
    We define the space $BMO\left(\Omega\right)=\left\{f\in L_\text{loc}^1\left(\Omega\right)\mid \left\lVert f\right\rVert_{BMO}<\infty\right\}$ together with the $BMO$ semi-norm as above.
\end{definition}
\begin{lemma}\cite[Lemma 1]{brenir2}
\label{bmol1}
Given $K\subset\Omega$ compact, there is a constant $C_K$ such that for every $f\in BMO\left(\Omega\right)$
\begin{equation}
    \left\lVert f-\Bar{f}_K\right\rVert_{L^1\left(K\right)}\leq C_K\left\lVert f\right\rVert_{BMO}
\end{equation}
where $\Bar{f}_K=\fint\limits_Kf\left(y\right)dy$.\\
Moreover, if $K$ is a closed ball $\overline{B}\subset\Omega$ then $C_{\overline{B}}=\lambda_n\left(B\right)$.
\end{lemma}
\begin{claim}
   $BMO\left(\Omega\right)$ modulo constant functions is a Banach Space.
\end{claim}
\begin{proof}\cite[67]{banach}
    It suffices to show that the $BMO$ semi-norm is complete. Let $\left\{f_n\right\}_{n=1}^\infty\subset BMO$ be a Cauchy sequence with respect to $BMO$ norm. 
    Given $K\subset\Omega$ compact, by Lemma \ref{bmol1} we have:
    $$\int\limits_K \left\lvert f_n\left(y\right)-\Bar{f}_{n,K}-\left(f_m\left(y\right)-\Bar{f}_{m,K}\right)\right\rvert dy\leq C_K\left\lVert f_n-f_m\right\rVert_{BMO}.$$
    Thus the sequence $\left\{f_n-\Bar{f}_{n,K}\right\}_{n=1}^\infty$ is Cauchy in $L^1\left(K\right)$, which gives the convergence 
    $$f_n-\bar{f}_{n,K}\overset{L^1\left(K\right)}{\longrightarrow}\Tilde{f}_K.$$ 
    Similarly, take some other $K'\subset\Omega$ compact, with $K\subset K'$, then we have 
    $$f_n-\bar{f}_{n,K'}\overset{L^1\left(K'\right)}{\longrightarrow}\Tilde{f}_{K'}.$$
    The last convergence also holds in $L^1\left(K\right)$, thus 
    $$\Bar{f}_{n,K}-\Bar{f}_{n,K'}\overset{L^1\left(K\right)}{\longrightarrow} C.$$
    But as constant functions this implies that as a sequence of numbers:
    $$\Bar{f}_{n,K}-\Bar{f}_{n,K'}\rightarrow C\left(K,K'\right).$$
    Now we define the limit function $f$ in the following way: \\
    Consider an ascending sequence of compact sets $\left\{K_i\right\}_{i=1}^\infty$ with $\bigcup\limits_{i=1}^\infty K_i=\Omega$.\\
    Given $x\in\Omega$, choose some $i\in\mathbb{N}$ with $x\in K_i$ and define $f\left(x\right)=\Tilde{f}_{K_i}-C(K_1,K_i)$. This is well defined: If we look at $j\in\mathbb{N}$ with $x\in K_j$ then we need $\Tilde{f}_{K_i}-C(K_1,K_i)=\Tilde{f}_{K_j}-C(K_1,K_j)$, assuming $j>i$, we need to show that:
    $$\Tilde{f}_{K_j}-\Tilde{f}_{K_i}=C(K_1,K_j)-C(K_1,K_i),$$
    which indeed follows from the above definitions.\\
    Now given some ball $B$ with $\overline{B}\subset\Omega$, take some $i\in\Omega$ with $B\subset K_i$, and we have that:
    \begin{equation*}
        \begin{split}
            &\int\limits_B \left\lvert f_n\left(y\right)-f\left(y\right)-\left(\Bar{f}_{n,B}-\Bar{f}_B\right)\right\rvert dy\\
            &\qquad =\int\limits_B \left\lvert f_n\left(y\right)-\Tilde{f}_{K_i}\left(y\right)+C\left(K_1,K_i\right)-\Bar{f}_{n,B}+\Bar{f}_B\right\rvert dy\\
            &\qquad = \int\limits_B \left\lvert f_n\left(y\right)-\Bar{f}_{n,B}-\Tilde{f}_B\left(y\right)+\left(\Tilde{f}_B\left(y\right)-\Tilde{f}_{K_i}\left(y\right)+C\left(K_1,K_i\right)+\Bar{f}_B\right)\right\rvert dy\\
            &\qquad = \int\limits_B \left\lvert f_n\left(y\right)-\Bar{f}_{n,B}-\Tilde{f}_B\left(y\right)+\left(-C\left(B,K_i\right)+C\left(K_1,K_i\right)+\Bar{f}_B\right)\right\rvert dy\\
             &\qquad = \int\limits_B \left\lvert f_n\left(y\right)-\Bar{f}_{n,B}-\Tilde{f}_B\left(y\right)+\left(-C\left(B,K_i\right)+{\overline{{\Tilde{f}}_{K_i}}}_B\right)\right\rvert dy
        \end{split}
    \end{equation*}
    We note that by definitions:
    $${\overline{\Tilde{f}_{K_i}}}_B=\fint\limits_B \Tilde{f}_{K_i}\left(y\right)dy=\underset{n\rightarrow\infty}{\lim}\fint\limits_B f_n\left(y\right)-\bar{f}_{n,K_i}=\underset{n\rightarrow\infty}{\lim}\bar{f}_{n,B}-\bar{f}_{n,K_i}=C(B,K_i)$$
    Thus:
    $$\int\limits_B \left\lvert f_n\left(y\right)-f\left(y\right)-\left(\Bar{f}_{n,B}-\Bar{f}_B\right)\right\rvert dy=\int\limits_B \left\lvert f_n\left(y\right)-\Bar{f}_{n,B}-\Tilde{f}_B\left(y\right)\right\rvert dy$$
    and since $f_m-\bar{f}_{m,B}\overset{L^1\left(B\right)}{\longrightarrow}\Tilde{f}_B$ we have that the right integral divided by the volume of the ball tends to zero uniformly on all the balls as:
      $$\fint\limits_B \left\lvert f_n\left(y\right)-\Bar{f}_{n,B}-\left(f_m\left(y\right)-\Bar{f}_{m,B}\right)\right\rvert dy\leq\left\lVert f_n-f_m\right\rVert_{BMO}$$
    This gives that $\left\lVert f_n-f\right\rVert_{BMO}\rightarrow0$ and that $f\in BMO$.
\end{proof}
\begin{remark}\cite[Remark 1]{brenir2}
    The space is independent of changing $r_x$ to
    $$r_x=\min\left\{k\cdot\dist\left(x,\Gamma\right),r_0\right\}\text{ where }0<k<1,\text{ and }r_0>0.$$ 
    Any such alternative choice gives an equivalent norm for the space. Additionally taking the averages over different shapes, such as cubes, also gives an equivalent norm.
\end{remark}
\begin{remark}\cite[Remark 2]{brenir2}
    Following the definitions, we obtain the continuous embedding 
    \begin{equation}
        \label{bmoemb}
        \left\lVert f\right\rVert_{BMO}\leq2\left\lVert f\right\rVert_{L^\infty}
    \end{equation}
\end{remark}
Using the above definition, we can define the $VMO$ space
\begin{definition}[$VMO$ Space]
    We define the space $VMO\left(\Omega\right)$ as:
    $$VMO\left(\Omega\right)=\overline{C\left(\overline{\Omega}\right)},$$
    where the closure is taken with respect to the $BMO$ norm.
\end{definition}
\begin{lemma}[Alternative Characterizations of $VMO$]\cite[Theorem 1]{brenir2}
\label{altvmo}
    The following statements are equivalent:
    \begin{enumerate}
        \item $f\in VMO\left(\Omega\right)$.
        \item $f\in BMO\left(\Omega\right)$, and $\fint\limits_{B_\epsilon\left(x\right)}\left\lvert f\left(y\right)-\Bar{f}_\epsilon\left(x\right)\right\rvert dy\underset{\epsilon\rightarrow0}{\longrightarrow}0$ uniformly on $x$.
        \item There is a sequence $\left\{f_i\right\}_{i=0}^\infty\subset C^\infty_0\left(\Omega\right)$ with $f_i\overset{BMO}{\longrightarrow}f$.
    \end{enumerate}
\end{lemma} 
With the above lemma we can prove the following:
\begin{claim}
\label{embedding}
    \cite[Example 1]{brenir2} 
    $$W^{1,n}\left(\Omega\right)\subset VMO\left(\Omega\right).$$
\end{claim}
\begin{claim}
    Given $f\in W^{1,n}\left(\Omega;\mathbb{R}^{n+1}\right)$ with $\nabla f$ of full rank, such that $\nu_f\in W^{1,n}\left(\Omega;S^n\right)$, then the function $\Tilde{F}\left(x,t\right):\Omega\times\left(-d,d\right)\rightarrow\mathbb{R}^{n+1}$ (for some small $1>d>0$) given by:
    $$\Tilde{F}\left(x,t\right)=f(x)+t\nu_f(x)$$ 
    satisfies $\Tilde{F}\in VMO\left(\Omega\times\left(-d,d\right);\mathbb{R}^{n+1}\right)\cap W^{1,n}\left(\Omega\times\left(-d,d\right);\mathbb{R}^{n+1}\right)$, and also $\det\nabla \Tilde{F}\in L^1\left(\Omega\times\left(-d,d\right);\mathbb{R}\right)$.
\end{claim}  
\begin{proof} 
    Given $\Tilde{x}=\left(x,t_0\right)$, and $\epsilon<\min\left\{1,\dist\left(\Tilde{x},\Gamma\times\left(-d,d\right)\right)\right\}$ by Poincare's Inequality there is a constant $C$ such that:
    $$\int\limits_{B_\epsilon\left(\Tilde{x}\right)}\left\lvert \Tilde{F}\left(\Tilde{y}\right)-\Bar{\Tilde{F}}_\epsilon\left(\Tilde{x}\right)\right\rvert d\Tilde{y}\leq C\epsilon\int\limits_{B_\epsilon\left(\Tilde{x}\right)}\left\lvert \nabla \Tilde{F}\left(\Tilde{y}\right)\right\rvert d\Tilde{y}.$$
    If we also note that $\nabla\Tilde{F}=\left(\nabla f+t\nabla\nu_f\mid\nu_f\right)$. Therefore we have that
    $$\int\limits_{B_\epsilon\left(\Tilde{x}\right)}\left\lvert \Tilde{F}\left(\Tilde{y}\right)-\Bar{\Tilde{F}}_\epsilon\left(\Tilde{x}\right)\right\rvert d\Tilde{y}\leq C\epsilon\int\limits_{B_\epsilon\left(\Tilde{x}\right)}\left\lvert \nabla f\right\rvert+\left\lvert t\nabla\nu_f\right\rvert+\left\lvert\nu_f\right\rvert d\Tilde{y}$$
    Now using the Fubini-Tonelli Theorem we have that
    $$\int\limits_{B_\epsilon\left(\Tilde{x}\right)}\left\lvert \Tilde{F}\left(\Tilde{y}\right)-\Bar{\Tilde{F}}_\epsilon\left(\Tilde{x}\right)\right\rvert d\Tilde{y}\leq C\epsilon\int\limits_{B_\epsilon\left(x\right)}\int\limits_{t_0-\sqrt{\epsilon^2-\left\lvert y\right\rvert^2}}^{t_0+\sqrt{\epsilon^2-\left\lvert y\right\rvert^2}}\left\lvert \nabla f\right\rvert+\left\lvert t\nabla\nu_f\right\rvert+\left\lvert\nu_f\right\rvert dtdy$$
    We can notice that $\sqrt{\epsilon^2-\left\lvert y\right\rvert^2}\leq\epsilon$, and also that
    $$\int\limits_{t_0-\sqrt{\epsilon^2-\left\lvert y\right\rvert^2}}^{t_0+\sqrt{\epsilon^2-\left\lvert y\right\rvert^2}}\left\lvert t\right\rvert dt=\left\{ \begin{array}{rcl}
         2\left\lvert t_0\right\rvert\sqrt{\epsilon^2-\left\lvert y\right\rvert^2} & \mbox{if}
         & \left\lvert t_0\right\rvert\geq\sqrt{\epsilon^2-\left\lvert y\right\rvert^2} \\ 
         t_0^2+\epsilon^2-\left\lvert y\right\rvert^2  & \mbox{if} & \left\lvert t_0\right\rvert<\sqrt{\epsilon^2-\left\lvert y\right\rvert^2}
                \end{array}\right.
                \leq 2\epsilon$$
    Note that in the above inequality, we used that $\left\lvert t_0\right\rvert<1$, and $\epsilon^2<\epsilon<1$. Thus we have that:
    $$\int\limits_{B_\epsilon\left(\Tilde{x}\right)}\left\lvert \Tilde{F}\left(\Tilde{y}\right)-\Bar{\Tilde{F}}_\epsilon\left(\Tilde{x}\right)\right\rvert d\Tilde{y}\leq C\epsilon^2\int\limits_{B_\epsilon\left(x\right)}\left\lvert \nabla f\right\rvert+\left\lvert \nabla\nu_f\right\rvert+\left\lvert\nu_f\right\rvert dy$$
    Finally, if we use the fact that $\left\lvert \nabla f\right\rvert+\left\lvert \nabla\nu_f\right\rvert+\left\lvert\nu_f\right\rvert\in L^n\left(\Omega\right)$ (as each individual function is) together with Hölder's Inequality we have that:
    \begin{equation*}
        \begin{split}
            &\int\limits_{B_\epsilon\left(\Tilde{x}\right)}\left\lvert \Tilde{F}\left(\Tilde{y}\right)-\Bar{\Tilde{F}}_\epsilon\left(\Tilde{x}\right)\right\rvert d\Tilde{y}\\
            &\qquad \leq C\epsilon^2\left(\lambda_n\left(B_1\left(0\right)\right)\epsilon^n\right)^\frac{n-1}{n}\left(\int\limits_{B_\epsilon\left(x\right)}\left(\left\lvert \nabla f\right\rvert+\left\lvert \nabla\nu_f\right\rvert+\left\lvert\nu_f\right\rvert\right)^n dy\right)^\frac{1}{n}.
        \end{split}
    \end{equation*}
    Thus we have that $\Tilde{F}\in BMO\left(\Omega\times\left(-d,d\right);\mathbb{R}^{n+1}\right)$ with bound
    $$\fint\limits_{B_\epsilon\left(\Tilde{x}\right)}\left\lvert \Tilde{F}\left(\Tilde{y}\right)-\Bar{\Tilde{F}}_\epsilon\left(\Tilde{x}\right)\right\rvert d\Tilde{y}\leq C\left(\int\limits_{B_\epsilon\left(x\right)}\left(\left\lvert \nabla f\right\rvert+\left\lvert \nabla\nu_f\right\rvert+\left\lvert\nu_f\right\rvert\right)^n dy\right)^\frac{1}{n}.$$
    By the uniformly integrable property of $L^1$ functions the right hand side is arbitrarily small for small $\epsilon$, independent of $x$. Using condition 2 of Lemma \ref{altvmo} we have that $\Tilde{F}\in VMO\left(\Omega\times\left(-d,d\right);\mathbb{R}^{n+1}\right)$.
    The other parts of the claim follow directly from the definitions.
\end{proof} 
\begin{remark}
    The argument for proving claim \ref{embedding} is very similar.
\end{remark}
Finally we present the following comparison lemma for averages of $VMO$ functions:
\begin{lemma}[Lemma A-B]\cite[Lemma A.4]{brenir1} 
\label{AB}
Let $\left(X,\mu\right)$, and $Y$ be measure spaces, then for any measurable function $g:X\rightarrow Y$ and measurable sets $A\subset B\subset X$ we have 
$$\left\lvert\Bar{g}_A-\Bar{g}_B\right\rvert\leq\frac{\mu\left(B\right)}{\mu\left(A\right)}\fint\limits_B\left\lvert g-\bar{g}_B\right\rvert d\mu .$$
\end{lemma}
\begin{proof}
    It follows from applying the integral triangle inequality:
    \begin{equation*}
        \begin{split}
            \left\lvert\Bar{g}_A-\Bar{g}_B\right\rvert &=\left\lvert\fint\limits_A gd\mu-\fint\limits_B gd\mu\right\rvert=\left\lvert\fint\limits_A gd\mu-\fint\limits_A\left(\fint\limits_B gd\mu\right)d\mu\right\rvert\\
            &=\left\lvert\fint\limits_A \left(g-\fint\limits_B gd\mu\right)d\mu\right\rvert\leq\fint\limits_A \left\lvert g-\fint\limits_B gd\mu\right\rvert d\mu\\
            &=\frac{\mu\left(B\right)}{\mu\left(B\right)}\fint\limits_A\left\lvert g-\bar{g}_B\right\rvert d\mu=\frac{\mu\left(B\right)}{\mu\left(B\right)\mu\left(A\right)}\int\limits_A\left\lvert g-\bar{g}_B\right\rvert d\mu\\
            &\leq\frac{\mu\left(B\right)}{\mu\left(A\right)\mu\left(B\right)}\int\limits_B\left\lvert g-\bar{g}_B\right\rvert d\mu=\frac{\mu\left(B\right)}{\mu\left(A\right)}\fint\limits_B\left\lvert g-\bar{g}_B\right\rvert d\mu
        \end{split}
    \end{equation*}
    as needed.
\end{proof}
\begin{corollary}
\label{avgc}
Let $f\in VMO\left(\Omega\right)$, then for every $\delta>0$, there is some $\epsilon_\delta>0$ such that for every $\epsilon_\delta>\epsilon>0$, and $x,y\in\Omega$ with $\dist\left(x,\Gamma\right)\geq3\epsilon$, and $\left\lvert x-y\right\rvert<2\epsilon$ we have:
$$\left\lvert \Bar{f}_\epsilon\left(x\right)-\Bar{f}_\epsilon\left(y\right)\right\rvert<\delta$$
\end{corollary}
\begin{proof}
    By Lemma \ref{altvmo} there is some $\epsilon_\delta>0$ such that for all $\epsilon_\delta>\epsilon>0$, and $x\in\Omega$ with $\dist\left(x,\Gamma\right)\geq3\epsilon$ we have 
    $$\fint\limits_{B_{3\epsilon}\left(x\right)}\left\lvert f\left(z\right)-\Bar{f}_{3\epsilon}\left(x\right)\right\rvert dz<\frac{\delta}{2\cdot3^n}.$$
    Using that $\left\lvert x-y\right\rvert<2\epsilon$ implies that $B_\epsilon\left(x\right),B_\epsilon\left(y\right)\subset B_{3\epsilon}\left(x\right)$, and applying the the triangle inequality, and Lemma \ref{AB} we get:
\begin{equation*}
    \begin{split}
        \left\lvert \Bar{f}_\epsilon\left(x\right)-\Bar{f}_\epsilon\left(y\right)\right\rvert&\leq\left\lvert \Bar{f}_\epsilon\left(x\right)-\Bar{f}_{3\epsilon}\left(x\right)\right\rvert+\left\lvert \Bar{f}_{3\epsilon}\left(x\right)-\Bar{f}_\epsilon\left(y\right)\right\rvert\\
        &\leq \frac{\lambda_n\left(B_{3\epsilon}\left(x\right)\right)}{\lambda_n\left(B_\epsilon\left(x\right)\right)}\fint\limits_{B_{3\epsilon}\left(x\right)}\left\lvert f\left(z\right)-\Bar{f}_{3\epsilon}\left(x\right)\right\rvert dz \\
        &+\frac{\lambda_n\left(B_{3\epsilon}\left(x\right)\right)}{\lambda_n\left(B_\epsilon\left(y\right)\right)}\fint\limits_{B_{3\epsilon}\left(x\right)}\left\lvert f\left(z\right)-\Bar{f}_{3\epsilon}\left(x\right)\right\rvert dz\\
        &=2\cdot3^n\fint\limits_{B_{3\epsilon}\left(x\right)}\left\lvert f\left(z\right)-\Bar{f}_{3\epsilon}\left(x\right)\right\rvert dz<2\cdot3^n\frac{\delta}{2\cdot3^n}=\delta
    \end{split}
\end{equation*}
    as needed.
\end{proof}
\subsection{Traces in $VMO$}
\label{tvmo}
By Claim \ref{embedding}, we have that $VMO$ space is larger than a critical Sobolev space. A natural step towards defining a consistent degree theory is to construct a notion of trace. Unfortunately, there is no consistent definition of trace for the whole $VMO\left(\Omega\right)$ space, but we can look at a smaller class of functions which do have a defined trace which is also in $VMO\left(\Gamma\right)$. First we define a basic case:
\begin{definition}[${VMO}_0$]
We say that $f\in VMO\left(\Omega\right)$ also has $f\in {VMO}_0\left(\Omega\right)$ if for any ball $B$ with $\Omega\subset B$ we have that:
 $$f_0\left(x\right)=\left\{ \begin{array}{rcl}
         f\left(x\right) & \mbox{if}
         & x\in\Omega \\ 
         0  & \mbox{if} & x\in B\backslash\Omega
                \end{array}\right.\in VMO\left(B\right).$$
\end{definition}
We have an alternative characterization for ${VMO}_0$. To write down this characterization we use the notation $D_\epsilon=\left\{x\in\Omega\mid\dist\left(x,\Gamma\right)=\epsilon\right\}$.
\begin{claim}\cite[Theorem 2]{brenir2}
\label{equiv}
    The following statements are equivalent:
    \begin{enumerate}
        \item $f\in {VMO}_0\left(\Omega\right).$
        \item $\lim\limits_{\epsilon\rightarrow0}\underset{x\in D_\epsilon}{\sup}\fint\limits_{B_{k\epsilon}\left(x\right)}\left\lvert f\left(y\right)\right\rvert dy=0$ for any $0<k\leq1$.
    \end{enumerate}
\end{claim}
\begin{remark}
    The above claim is only true when there are some $\epsilon_0,\alpha>0$ such that for every $y\in\Gamma$, and $0<\epsilon<\epsilon_0$ we have:   $$\lambda_n\left(B_{\epsilon}\left(y\right)\cap\Omega^C\right)\geq\alpha\lambda_n\left(B_{\epsilon}\left(y\right)\right).$$
    This is clearly satisfied if $\Gamma$ is Lipschitz or smooth.
\end{remark}
To define the more general ${VMO}_\varphi\left(\Omega\right)$ for $\varphi\in VMO\left(\Gamma\right)$ we need the following lemma:
\begin{lemma} \cite[Lemma 5]{brenir2}
\label{ext}
    Given $\varphi\in VMO\left(\Gamma
    \right)$, then the function $\Bar{\varphi}=\varphi\circ P$ is in $VMO\left(T\right)$ (recall that $T$ is a tubular neighborhood of $\Gamma$, and $P$ is the projection operator from $T$ to $\Gamma$).
\end{lemma}
\begin{remark}
On \cite{brenir2} the proof of this is just the first part of the proof of the Lemma stated there.
\end{remark}
Let $\psi:\Omega\rightarrow\mathbb{R}$ be a smooth function, such that $\psi=0$ outside of $T$, and $\psi=1$ for some open $U\subset\Omega$ with $\Gamma\subset\overline{U}\subset T$. Then it follows from the lemma above that $\Tilde{\varphi}=\psi\cdot\Bar{\varphi}\in VMO\left(\Omega\right)$.
\begin{definition}[${VMO}_\varphi$]
We say that $f\in VMO\left(\Omega\right)$ also satisfies $f\in {VMO}_\varphi\left(\Omega\right)$ if we have that:
$$\Bar{f}_\varphi\left(x\right)=f\left(x\right)-\Tilde{\varphi}\left(x\right)\in {VMO}_0\left(\Omega\right).$$
\end{definition}
\begin{remark}
    An equivalent definition is the following:
    $$f\in {VMO}_\varphi\left(\Omega\right)\iff f_\varphi\left(x\right)=\left\{ \begin{array}{rcl}
         f\left(x\right) & \mbox{if}
         & x\in\Omega \\ 
        \Bar{\varphi}\left(x\right)  & \mbox{if} & x\in T\backslash\Omega
                \end{array}\right.\in VMO\left(\Omega\cup T\right)$$
                Where $T$ is as in Lemma \ref{ext}.
\end{remark}
\begin{claim}
\label{codim}
    Let $f\in W^{1,n}\left(\Omega;\mathbb{R}^{n+1}\right)$ such that $\trace_{\Gamma}f\in W^{1,n}\left(\Gamma;\mathbb{R}^{n+1}\right)$ ,and $\nabla f$ is of full rank, with $\nu_f\in W^{1,n}\left(\Omega;S^n\right)$, and $\trace_\Gamma\nu_f\in W^{1,n}\left(\Gamma;S^n\right)$. Then the function $\Tilde{F}\left(x,t\right):\Omega\times\left(-d,d\right)\rightarrow\mathbb{R}^{n+1}$ (for some small $1>d>0$) given by:
    $$\Tilde{F}\left(x,t\right)=f(x)+t\nu_f(x)$$ 
    satisfies $\Tilde{F}\in {VMO}_{\varphi_f}\left(\Omega\times\left(-d,d\right);\mathbb{R}^{n+1}\right)$ where $\varphi_f$ id given by:
    $$\varphi_f\left(x,t\right)=\left\{ \begin{array}{rcl}
         f\left(x\right)+d\nu_f\left(x\right) & \mbox{if}
         & \left(x,t\right)\in\Omega\times\left\{d\right\} \\ 
         f\left(x\right)-d\nu_f\left(x\right)  & \mbox{if} & \left(x,t\right)\in\Omega\times\left\{-d\right\} \\
         \trace_\Gamma f\left(x\right)+t\trace_\Gamma\nu_f\left(x\right)  & \mbox{if} & \left(x,t\right)\in \Gamma\times\left[-d,d\right]
                \end{array}\right..$$
\end{claim}
\begin{proof}
    It suffices to show that $\varphi_f\in VMO\left(\partial\left(\Omega\times\left(-d,d\right)\right)\right)$ which is clear from claim \ref{embedding}.
\end{proof}
In order to define the meaning of an element being outside the trace, we need the concept of essential range:
\begin{definition}[Essential Range]
    Given a measure space $\left(X,\mu\right)$, a topological space $\left(Y,\tau\right)$ (with Borel measure), and a measurable map $f:X\rightarrow Y$, the essential range of the map is given by:
    $$\eim\left(f\right)=\left\{y\in Y\mid0<\mu\left(f^{-1}\left(S\right)\right) \text{ for all }y\in S\in\tau\right\}.$$
\end{definition}
\begin{lemma}
\label{bound1}
    Let $f\in{VMO}_\varphi\left(\Omega\right)$, and $p\in\mathbb{R}^n\backslash\overline{\eim\left(\varphi\right)}$, then there is some $d_0>0$ and a neighborhood $N$ of $\Gamma$ in $\Omega$ such that $\fint\limits_B\left\lvert f\left(y\right)-p\right\rvert dy> d_0$ for every $$B\in\left\{B=B_\epsilon\left(x\right)\subset N\mid \epsilon=\dist\left(x,\Gamma\right)\right\}=:D_N.$$ 
\end{lemma}
\begin{proof}
    Without loss of generality, we may assume that $p=0$. Let $2d_0=\dist\left(0,\eim\left(\varphi\right)\right)>0$ thus, for a.e.\ $x\in\Gamma$ we have that $\left\lvert\varphi\left(x\right)\right\rvert\geq2d_0$. From the assumption that $f\in {VMO}_\varphi\left(\Omega\right)$, and Claim \ref{equiv}, we have some $\epsilon_0>0$ such that for all $0<\epsilon<\epsilon_0$, $\underset{x\in D_\epsilon}{\sup}\fint_{B_\epsilon\left(x\right)}\left\lvert\Bar{f}_\varphi\left(y\right)\right\rvert dy< d_0$. \\
    Recall that we defined $U$ to be the neighborhood of $\Gamma$ in $\Omega$ with $\Tilde{\varphi}=\Bar{\varphi}$, and let $N=\left\{x\in\Omega\mid B_\epsilon\left(x\right)\subset U\text{ where } \epsilon=\dist\left(x,\Gamma\right)<\epsilon_0\right\}$. If $B\in D_N$, then by the inverse triangle inequality we have:
    \begin{equation*}
        \begin{split}
            d_0 & >\fint_B\left\lvert\Bar{f}_\varphi\left(y\right)\right\rvert dy=\fint_B\left\lvert f\left(y\right)-\Tilde{\varphi}\left(y\right)\right\rvert dy=\fint_B\left\lvert f\left(y\right)-\Bar{\varphi}\left(y\right)\right\rvert dy \\
            & \geq\fint_B\left\lvert\Bar{
            \varphi}\left(y\right)\right\rvert dy-\fint_B\left\lvert f\left(y\right)\right\rvert dy\geq 2d_0-\fint_B\left\lvert f\left(y\right)\right\rvert dy,
        \end{split}
\end{equation*}
which gives the desired inequality: $\fint_B\left\lvert f\left(y\right)\right\rvert dy> d_0$ for every $B\in D_N$.
\end{proof}
\subsection{Approximation by smooth functions}
In this subsection we want to construct approximations of $VMO$ functions by smooth ones. This approximations will subsequently allow us to properly define a degree in $VMO$, and prove properties of it.\\
The first approximation we have is:
\begin{claim}\cite[Corollary 1]{brenir1}
\label{avgconv}
    Let $f\in VMO\left(\Omega\right)$ then:
    $$\lim\limits_{\epsilon\rightarrow0}\left\lVert f-\Bar{f}_\epsilon\right\rVert_{BMO}=0.$$
\end{claim}
To improve this, we need the following lemma about $VMO$ functions:
\begin{lemma}
\label{nullker}
Let $\psi\in C^\infty\left(\overline{B_1\left(0\right)}\right)$ such that     $\int\limits_{B_1\left(0\right)}\psi\left(z\right)dz=0$, and let $f\in VMO\left(\Omega\right)$, then:
$$\fint\limits_{B_\epsilon\left(x\right)}f\left(z\right)\psi\left(\frac{x-z}{\epsilon}\right)dz\underset{\epsilon\rightarrow0}{\longrightarrow}0$$
uniformly on $x\in\Omega$.
\end{lemma}
\begin{proof}
    Let $\delta>0$, as $f\in VMO\left(\Omega\right)$, there is some $\epsilon_0$ such that for every $\epsilon<\epsilon_0$, and $x\in\Omega$ with $B_\epsilon\left(x\right)\subset\Omega$ we have:
    $$\fint\limits_{B_\epsilon\left(x\right)}\left\lvert f\left(y\right)-\Bar{f}_\epsilon\left(x\right)\right\rvert dy<\frac{\delta}{\left\lVert\psi\right\rVert_\infty}.$$
    Then by using $\int\limits_{B_1\left(0\right)}\psi\left(z\right)dz=0$, the integral triangle inequality, and Holder's inequality we get:
    \begin{equation*}
\begin{split}
&\left\lvert\fint\limits_{B_\epsilon\left(x\right)}f\left(z\right)\psi\left(\frac{x-z}{\epsilon}\right)dz\right\rvert \\
&\qquad=\left\lvert\fint\limits_{B_\epsilon\left(x\right)}f\left(z\right)\psi\left(\frac{x-z}{\epsilon}\right)dz-\fint\limits_{B_\epsilon\left(x\right)}\Bar{f}_\epsilon\left(x\right)\psi\left(\frac{x-z}{\epsilon}\right)dz\right\rvert\\
&\qquad = \left\lvert\fint\limits_{B_\epsilon\left(x\right)}\left(f\left(z\right)-\Bar{f}_\epsilon\left(x\right)\right)\psi\left(\frac{x-z}{\epsilon}\right)dz\right\rvert\\
&\qquad \leq \fint\limits_{B_\epsilon\left(x\right)}\left\lvert\left(f\left(z\right)-\Bar{f}_\epsilon\left(x\right)\right)\psi\left(\frac{x-z}{\epsilon}\right)\right\rvert dz\\
&\qquad \leq\left\lVert\psi\right\rVert_\infty\fint\limits_{B_\epsilon\left(x\right)}\left\lvert f\left(y\right)-\Bar{f}_\epsilon\left(x\right)\right\rvert dy<\left\lVert\psi\right\rVert_\infty\frac{\delta}{\left\lVert\psi\right\rVert_\infty}=\delta.
        \end{split}
    \end{equation*}
Because $\epsilon$ is chosen uniformly for all $x\in\Omega$, the proof is complete.
\end{proof}
For the last theorem of the section we will use mollifiers, hence we define a mollificating function.
\begin{definition}
    We define $\eta:B_1\left(0\right)\rightarrow\mathbb{R}_{\geq0}$ to be a function satisfying:
    \begin{enumerate}
        \item $\eta\in C_c^\infty\left(B_1\left(0\right)\right)$.
        \item $\int\limits_{B_1\left(0\right)}\eta\left(y\right)dy=1$.
        \item $\eta$ is radially symmetric.
    \end{enumerate}
    Additionally, for $\epsilon>0$, we define $\eta_\epsilon\left(x\right)=\frac{1}{\epsilon^n}\eta\left(\frac{x}{\epsilon}\right)\in C_c^\infty\left(B_\epsilon\left(0\right)\right)$.
\end{definition}
\begin{remark}
    Let $\omega_n=\lambda_n\left(B_1\left(0\right)\right)$, then the function $\psi_n=1-\omega_n\eta$ satisfies the same conditions as $\psi$ in Lemma \ref{nullker}.
\end{remark}
\begin{lemma}
\label{mollify}
    Recall that $T$ is a tubular neighborhood of $\Gamma$. If $g\in VMO\left(\Omega\cup T\right)$, then there is some $\epsilon_0>0$ such that for every $\epsilon_0>\epsilon>0$, $g*\eta_\epsilon$ is defined on $\Omega$, and
    $$\lim\limits_{\epsilon\rightarrow0}\left\lVert g*\eta_\epsilon-g\right\rVert_{BMO\left(\Omega\right)}=0.$$
\end{lemma}
\begin{proof}
    Because $T$ is a tubular neighborhood of $\Gamma$, we can find some $\epsilon_0>0$ such that for every $x\in\Omega$, and $0<\epsilon<\epsilon_0$ we have $B_\epsilon\left(x\right)\subset\Omega\cup T$ and thus $g*\eta_\epsilon\left(x\right)$ is defined.\\
    Let $\delta>0$, by Claim \ref{avgconv} there is some $\epsilon_1$ such that for every $0<\epsilon<\epsilon_1$ we have 
    $$\left\lVert g-\Bar{g}_\epsilon\right\rVert_{BMO}<\frac{\delta}{2}$$
    Let $x\in\Omega$, and let $0<\epsilon<\epsilon_0$, then we have that:
    \begin{equation}
    \label{avgmol}
    \begin{split}
        \Bar{g}_\epsilon\left(x\right)-g*\eta_\epsilon\left(x\right) & =\fint\limits_{B_\epsilon\left(x\right)}g\left(y\right)dy-\int\limits_{B_\epsilon\left(x\right)}g\left(y\right)\eta_\epsilon\left(x-y\right)dy \\
        & =\fint\limits_{B_\epsilon\left(x\right)}g\left(y\right)dy-\int\limits_{B_\epsilon\left(x\right)}g\left(y\right)\frac{1}{\epsilon^n}\eta\left(\frac{x-y}{\epsilon}\right)dy\\
        &=\fint\limits_{B_\epsilon\left(x\right)}g\left(y\right)dy-\int\limits_{B_\epsilon\left(x\right)}g\left(y\right)\frac{\omega_n}{\omega_n\epsilon^n}\eta\left(\frac{x-y}{\epsilon}\right)dy\\
        &=\fint\limits_{B_\epsilon\left(x\right)}g\left(y\right)dy-\frac{1}{\lambda_n\left(B_\epsilon\left(x\right)\right)}\int\limits_{B_\epsilon\left(x\right)}g\left(y\right)\omega_n\eta\left(\frac{x-y}{\epsilon}\right)dy\\
        &=\fint\limits_{B_\epsilon\left(x\right)}g\left(y\right)\left(1-\omega_n\eta\left(\frac{x-y}{\epsilon}\right)\right)dy\\ 
        &=\fint\limits_{B_\epsilon\left(x\right)}g\left(y\right)\psi_n\left(\frac{x-y}{\epsilon}\right)dy.
    \end{split}
    \end{equation}
    Thus, by Lemma \ref{nullker}, there is some $\epsilon_2>0$ such that for every $0<\epsilon<\epsilon_2$ we have:
    \begin{equation*}
        \begin{split}
            \left\lVert\Bar{g}_\epsilon\left(x\right)-g*\eta_\epsilon\left(x\right)\right\rVert_{L^\infty\left(\Omega\right)}&=\sup\limits_{x\in\Omega}\left\lvert\Bar{g}_\epsilon\left(x\right)-g*\eta_\epsilon\left(x\right)\right\rvert\\
            &=\sup\limits_{x\in\Omega}\left\lvert\fint\limits_{B_\epsilon\left(x\right)}g\left(y\right)\psi_n\left(\frac{x-y}{\epsilon}\right)dy\right\rvert<\frac{\delta}{4}.
        \end{split}
    \end{equation*}
    Therefore, for every $0<\epsilon<\min\left\{\epsilon_0,\epsilon_1,\epsilon_2\right\}$, by the triangle inequality, and Inequality \eqref{bmoemb} we have:
\begin{equation*}
    \begin{split}
        \left\lVert g*\eta_\epsilon-g\right\rVert_{BMO\left(\Omega\right)}&\leq\left\lVert g*\eta_\epsilon-\Bar{g}_\epsilon\right\rVert_{BMO\left(\Omega\right)}+\left\lVert \Bar{g}_\epsilon-g\right\rVert_{BMO\left(\Omega\right)}\\
        &\leq 2\left\lVert g*\eta_\epsilon-\Bar{g}_\epsilon\right\rVert_{L^\infty\left(\Omega\right)}+\left\lVert \Bar{g}_\epsilon-g\right\rVert_{BMO\left(\Omega\right)}\\
        &<2\frac{\delta}{4}+\frac{\delta}{2}=\delta
    \end{split}
\end{equation*}
    As needed.
\end{proof}
\begin{definition}[Extension by Reflection]
    Given $f:\Omega\rightarrow\mathbb{R}^n$ define $f_T:\Omega\cup T\rightarrow\mathbb{R}^n$ by:
    $$f_T\left(x\right)=\left\{ \begin{array}{rcl}
         f\left(x\right) & \mbox{if}
         & x\in\Omega \\ 
        f\left(2P\left(x\right)-x\right)  & \mbox{if} & x\in V\backslash\Omega
                \end{array}\right..$$
    (Recall that $P:T\rightarrow\Gamma$ is the projection operator.)
\end{definition}
\begin{remark}
    The way $T$ was defined on Section \ref{notation} gives that $f_T$ is well defined. 
\end{remark}
\begin{claim} 
    If $f\in VMO\left(\Omega\right)$ then $f_T\in VMO\left(\Omega\cup T\right)$.
\end{claim}
\begin{proof}
    Follows immediately from the definitions.
\end{proof}
\begin{theorem}
\label{bound2}
Given $f\in {VMO}_\varphi\left(\Omega\right)$, and $p\in\mathbb{R}^n\backslash\overline{\eim\left(\varphi\right)}$, there is some $c>0$ such that $f_i=f_T*\eta_{\frac{c}{i}}$ is well defined for all $i\in\mathbb{N}$, and there is some $d_1>0$, and a neighborhood $M$ of $\Gamma$ in $\Omega$ with $\fint\limits_B\left\lvert f_i\left(y\right)-p\right\rvert dy> d_1$ for every $B\in D_M$, and $i\in\mathbb{N}$.
\end{theorem}
\begin{proof}
Let $d_0>0$, $N$ as obtained from Lemma \ref{bound1}, and define $d_1=\frac{d_0}{4}$. Define some $c_0>0$ such that for every $x\in\Omega$ and $c_0>\epsilon>0$ we have $B_{3\epsilon}\left(x\right)\subset\Omega\cup T$.\\
By Lemma \ref{altvmo} there is some $c_1>0$ such that for every $x\in\Omega$, and $\epsilon<c_1$ with $B_\epsilon\left(x\right)\subset\Omega$ we have
    $$\fint\limits_{B_\epsilon\left(x\right)}\left\lvert f\left(y\right)-\Bar{f}_\epsilon\left(x\right)\right\rvert dy<d_1.$$
    We also define $c_1$ such that $\left\{x\in\Omega\mid\dist\left(x,\Gamma\right)\leq 2c_1\right\}\subset N$. If we let $\epsilon<c_1$, and $x\in\Omega$ with $\dist\left(x,\Gamma\right)=\epsilon$ then $B_{\epsilon}\left(x\right)\in D_N$, and by the inverse triangle inequality, we have 
    \begin{equation*}
        \begin{split}
        \left\lvert\Bar{f}_\epsilon\left(x\right)-p\right\rvert=&\fint\limits_{B_\epsilon\left(x\right)}\left\lvert \Bar{f}_\epsilon\left(x\right)-p\right\rvert dy\\
        &\geq\fint\limits_{B_\epsilon\left(x\right)}\left\lvert f\left(y\right)-p\right\rvert dy-\fint\limits_{B_\epsilon\left(x\right)}\left\lvert f\left(y\right)-\Bar{f}_\epsilon\left(x\right)\right\rvert dy> 3d_1.
        \end{split}
    \end{equation*}
    Thus, as $B_\epsilon\left(x\right)\subset\Omega$ we have:
    $$\left\lvert\overline{f_T}_\epsilon\left(x\right)-p\right\rvert=\left\lvert\Bar{f}_\epsilon\left(x\right)-p\right\rvert>3d_1.$$
    By Equation \eqref{avgmol}, and Lemma \ref{nullker}, there is some $c_2>0$ such that for every $0<\epsilon<c_2$ and $x\in\Omega$ we have:
    $$\left\lvert\overline{f_T}_\epsilon\left(x\right)-f_T*\eta_\epsilon\left(x\right)\right\rvert<d_1.$$
    Using Corollary \ref{avgc}, there is some $c_3>0$ such that for every $c_3>\epsilon>0$, and $x,y\in\Omega$ with $B  _{3\epsilon}\left(x\right)\subset\Omega\cup T$, and $\left\lvert x-y\right\rvert<2\epsilon$ we have: 
    \begin{equation}
        \label{avgosc}
        \left\lvert\overline{f_T}_\epsilon\left(x\right)-\overline{f_T}_\epsilon\left(y\right)\right\rvert<d_1.
    \end{equation}
    By Lemma \ref{mollify} there is some $c_4>0$ such that for every $0<\epsilon<c_4$, and $B\subset\Omega$, we have:
    $$\fint\limits_B\left\lvert \left(f-f_T*\eta_\epsilon\right)\left(y\right)-\overline{f-f_T*\eta_\epsilon}_B\right\rvert dy\leq\left\lVert f-f_T*\eta_\epsilon\right\rVert_{BMO}<d_1.$$
    Choose any $0<c<\min\left\{c_0,c_1,c_2,c_3,c_4\right\}$, and define $M=\left\{x\in\Omega\mid\dist\left(x,\Gamma\right)\leq2c\right\}$ (thus $D_M\subset D_N$), and for every $i\in\mathbb{N}$, define $f_i=f_T*\eta_{\frac{c}{i}}$. For every $i\in\mathbb{N}$, we have the following:\\
    If $y\in\Omega$ with $\dist\left(y,\Gamma\right)\leq\frac{c}{i}$, there is a $x\in\Omega$ with $\dist\left(x,\Gamma\right)=\frac{c}{i}$ and $\left\lvert x-y\right\rvert<\frac{c}{i}<2\frac{c}{i}$, thus by the triangle inequality, and the above inequalities, we have:
    \begin{equation*}
        \begin{split}
        \left\lvert f_i\left(y\right)-p\right\rvert&\geq\left\lvert \overline{f_T}_\frac{c}{i}\left(x\right)-p\right\rvert-\left\lvert \overline{f_T}_\frac{c}{i}\left(x\right)-f_i\left(y\right)\right\rvert\\
        &\geq\left\lvert \overline{f_T}_\frac{c}{i}\left(x\right)-p\right\rvert-\left\lvert \overline{f_T}_\frac{c}{i}\left(x\right)-\overline{f_T}_\frac{c}{i}\left(y\right)\right\rvert-\left\lvert \overline{f_T}_\frac{c}{i}\left(y\right)-f_i\left(y\right)\right\rvert\\
        &>3d_1-d_1-d_1=d_1  
        \end{split}
    \end{equation*}
    Hence, if $B\in D _M$ with $B\subset\left\{x\in\Omega\mid\dist\left(x,\Gamma\right)\leq\frac{c}{i}\right\}$ then $\fint\limits_B\left\lvert f_i\left(y\right)-p\right\rvert dy> d_1$.\\
    If $B\in D_M$ with $B\not\subset\left\{x\in\Omega\mid\dist\left(x,\Gamma\right)\leq\frac{c}{i}\right\}$, then, there is some $x\in\Omega$ with $r=\dist\left(x,\Gamma\right)\leq c$ such that $B=B_r\left(x\right)$, and $\frac{c}{2i}<r$. By applying the triangle inequality we have:
    \begin{equation*}
        \begin{split}
            \fint\limits_B\left\lvert f_i\left(y\right)-p\right\rvert dy&\geq\fint\limits_B\left\lvert f\left(y\right)-p\right\rvert dy-\fint\limits_B\left\lvert f\left(y\right)-f_i\left(y\right)\right\rvert dy\\
            &\geq\fint\limits_B\left\lvert f\left(y\right)-p\right\rvert dy-\fint\limits_B\left\lvert \left(f-f_i\right)\left(y\right)-\overline{f-f_i}_B\right\rvert dy-\fint\limits_B\left\lvert \overline{f-f_i}_B\right\rvert dy\\
            &=\fint\limits_B\left\lvert f\left(y\right)-p\right\rvert dy-\fint\limits_B\left\lvert \left(f-f_i\right)\left(y\right)-\overline{f-f_i}_B\right\rvert dy-\left\lvert \overline{f-f_i}_B\right\rvert\\
            &>4d_1-d_1-d_1=2d_1>d_1
        \end{split}
    \end{equation*}
    as needed. In the last step, we used the bound $\left\lvert\overline{f-f_i}_B\right\rvert<d_1$, this follows from the fact that $\int\limits_{B_\frac{c}{i}\left(0\right)}\eta_\frac{c}{i}\left(y\right)dy=1$, the Fubini Tonelli theorem, the integral triangle inequality, and Inequality \eqref{avgosc} (as $\frac{c}{i}<2r\leq2c$):
    \begin{equation*}
        \begin{split}
            \left\lvert\overline{f-f_i}_B\right\rvert&=\left\lvert\overline{f_T-f_i}_B\right\rvert=\left\lvert\fint\limits_Bf_T\left(y\right)-f_i\left(y\right)dy\right\rvert\\
            &=\left\lvert\fint\limits_B\left(\int\limits_{B_\frac{c}{i}\left(0\right)}\left(f_T\left(y\right)-f_T\left(y-z\right)\right)\eta_\frac{c}{i}\left(z\right)dz\right)dy\right\rvert\\
            &=\left\lvert\int\limits_{B_\frac{c}{i}\left(0\right)}\left(\fint\limits_Bf_T\left(y\right)-f_T\left(y-z\right)dy\right)\eta_\frac{c}{i}\left(z\right)dz\right\rvert\\
            &\leq\int\limits_{B_\frac{c}{i}\left(0\right)}\left\lvert\fint\limits_Bf_T\left(y\right)-f_T\left(y-z\right)dy\right\rvert\eta_\frac{c}{i}\left(z\right)dz\\
            &=\int\limits_{B_\frac{c}{i}\left(0\right)}\left\lvert\fint\limits_{B_r\left(x\right)}f_T\left(y\right)-f_T\left(y-z\right)dy\right\rvert\eta_\frac{c}{i}\left(z\right)dz\\
            &=\int\limits_{B_\frac{c}{i}\left(0\right)}\left\lvert\overline{f_T}_r\left(x\right)-\overline{f_T}_r\left(x-z\right)\right\rvert\eta_\frac{c}{i}\left(z\right)dz\\
            &\leq\int\limits_{B_\frac{c}{i}\left(0\right)}d_1\eta_\frac{c}{i}\left(z\right)dz=d_1
        \end{split}
    \end{equation*}
    giving the desired bound.
\end{proof}
\begin{remark}
    The same proof works if we just assume that $f\in VMO\left(\Omega\right)$, and $p\in\mathbb{R}^n$ such that there are $d_0>0$, and $N$ neighborhood of $\Gamma$ in $\Omega$ such that $\fint\limits_B\left\lvert f\left(y\right)-p\right\rvert dy> d_0$ for every $B\in D_N$. 
\end{remark}
\subsection{Degree in $VMO$}
In Section \ref{tvmo} we gave a characterization of the boundary values for a map of class $VMO$, thus in this section we will define the degree for such maps. First we introduce the notation:
$$\Omega_\epsilon=\left\{x\in\Omega\mid \dist\left(x,\Gamma\right)>\epsilon\right\}.$$
With this notation, we can define the degree for $VMO$ maps.
\begin{definition}[$VMO$ Degree]
    Let $f\in VMO\left(\Omega\right)$, and let $p\in\mathbb{R}^n$ such that there is some $d_0>0$ and a neighborhood $N$ of $\Gamma$ in $\Omega$ such that $\fint\limits_B\left\lvert f\left(y\right)-p\right\rvert dy> d_0$ for every $B\in D_N$, then we can define the degree of $f$ at the point $p$ by:
$$\deg\left(f,\Omega,p\right)=\lim\limits_{\epsilon\rightarrow0}\deg\left(\Bar{f}_\epsilon,\Omega_{\epsilon},p\right),$$
where we recall that $\bar{f}_\epsilon\left(x\right)=\fint\limits_{B_\epsilon\left(x\right)}f\left(y\right)dy$.
\end{definition}
\begin{claim}\cite{brenir2}
    The degree for a $VMO$ map is well defined.
\end{claim}
\begin{proof}
   As done in the proof of Theorem \ref{bound2}, there is some $\epsilon_0>0$ such that for every $0<\epsilon<\epsilon_0$, and for every $x\in\partial\Omega_\epsilon$ we have $\left\lvert\Bar{f}_\epsilon\left(x\right)-p\right\rvert>\frac{d_0}{2}$, thus $\deg\left(\Bar{f}_\epsilon,\Omega_{\epsilon},p\right)$ is well defined, as $\Bar{f}_\epsilon$ is continuous. To prove the convergence of the limit, we just need to note that $\Bar{f_\epsilon}\left(x\right)$ is continuous where it is defined (as a function $\Omega\times\left(0,1\right)\rightarrow\mathbb{R}^n$), and use the fact that degree is homotopy invariant, and has an excision property.
\end{proof}
\begin{claim}[Properties of $VMO$ degree]\cite{brenir2}
\label{propvmodeg}
    Given $f\in VMO\left(\Omega\right)$, and $p\in\mathbb{R}^n$ as in the definition of degree, we have:
    \begin{enumerate}
        \item $\deg\left(f,\Omega,p\right)\neq0$ then $p\in\eim\left(f\right)$.
        \item If $f_i\rightarrow f$ both in $BMO$ and $L^1_\text{loc}$, and there exists $d_1>0$, and a neighborhood $M$ in $\Omega$ of $\Gamma$ with $\fint\limits_B\left\lvert f_i\left(y\right)-p\right\rvert dy>d_1$ for every $B\in D_M$, and $i\in\mathbb{N}$ then $\deg\left(f,\Omega,p\right)=\lim\limits_{i\rightarrow\infty}\deg\left(f_i,\Omega,p\right)$.
        \item Let $g\in VMO\left(\Omega\right)$ such that there is some neighborhood $\Tilde{N}$ of $\Gamma$ in $\Omega$, and some $d'<d_0$ with $\fint\limits_{B}\left\lvert f\left(y\right)-g\left(y\right)\right\rvert dy< d'$ for every $B\in D_{\Tilde{N}}$, then $\deg\left(f,\Omega,p\right)=\deg\left(g,\Omega,p\right)$.
        \item If $f\in {VMO}_\varphi\left(\Omega\right)$, $q\in\mathbb{R}^n$ is in the same connected component of $\mathbb{R}^n\backslash\overline{\eim\left(\varphi\right)}$ as $p$ then $\deg\left(f,\Omega,q\right)=\deg\left(f,\Omega,p\right)$.
    \end{enumerate}
\end{claim}
\begin{remark}
    The second property means that the degree is stable in the $VMO$ topology, and an equivalent version of homotopy invariance can be stated.\\
    The third property means that for any two maps $f,g\in{VMO}_\varphi\left(\Omega\right)$, and $p\in\mathbb{R}^n\backslash\overline{\eim\left(\varphi\right)}$ we have:
    $$\deg\left(f,\Omega,p\right)=\deg\left(g,\Omega,p\right).$$
\end{remark}
We recall Pratt's theorem:
\begin{theorem}[Pratt's Theorem]\cite[Theorem A.10]{rindler}
    Let $\left\{g_n\right\}_{n=1}^\infty$, be Lebesgue-measurable functions. If $g_n\rightarrow g$ either point-wise a.e.\ or in measure, and there is $\left\{h_n\right\}_{n=1}^\infty\subset L^1\left(\Omega\right)$ with $h_n\overset{L^1}{\longrightarrow}h$,and $\left\lvert g_n\right\rvert\leq h_n$ then $g_n\overset{L^1}{\longrightarrow}g$.
\end{theorem}
To prove the next theorem, we need the following claim and its corollary:
\begin{claim}
Let $\left\{g^1_n\right\}_{n=1}^\infty,...,\left\{g^k_n\right\}_{n=1}^\infty\subset L^p$ with $p\geq k$, and assume that $g^i_n\overset{L^p}{\rightarrow}g^i$ then $\prod\limits_{i=1}^kg^i_n\overset{L^\frac{p}{k}}{\rightarrow}\prod\limits_{i=1}^kg^i$.
\end{claim}
\begin{proof}
    Let's see the result by induction on $k.\\$
    The case $k=1$ is trivial. \\
    Assume that $\prod\limits_{i=1}^{k-1}g^i_n\overset{L^\frac{p}{k-1}}{\rightarrow}\prod\limits_{i=1}^{k-1}g^i$, and let $\delta>0$. By convergence there is some $M\in\mathbb{N}$ such that for every $n\in\mathbb{N}$ we have $\left\lVert g^k_n\right\rVert_p,\left\lVert \prod\limits_{i=1}^{k-1}g^i_n\right\rVert_\frac{p}{k-1}<M$, and some $N\in\mathbb{N}$ such that every $n>N$ has $\left\lVert g^k_n-g^k\right\rVert_p,\left\lVert \prod\limits_{i=1}^{k-1}g^i_n-\prod\limits_{i=1}^{k-1}g^i\right\rVert_\frac{p}{k-1}<\frac{\delta}{2M}$. By the triangle inequality, and Hölder's inequality we have:
    \begin{equation*}
    \begin{split}
        \left\lVert \prod\limits_{i=1}^kg^i_n-\prod\limits_{i=1}^kg^i\right\rVert_\frac{p}{k}&=\left\lVert g^k_n\prod\limits_{i=1}^{k-1}g^i_n-g^k\prod\limits_{i=1}^{k-1}g^i\right\rVert_\frac{p}{k}\\
        &\leq\left\lVert g^k_n\prod\limits_{i=1}^{k-1}g^i_n-g^k_n\prod\limits_{i=1}^{k-1}g^i\right\rVert_\frac{p}{k}+\left\lVert g^k_n\prod\limits_{i=1}^{k-1}g^i-g^k\prod\limits_{i=1}^{k-1}g^i\right\rVert_\frac{p}{k}\\
        &=\left\lVert g^k_n\left(\prod\limits_{i=1}^{k-1}g^i_n-\prod\limits_{i=1}^{k-1}g^i\right)\right\rVert_\frac{p}{k}+\left\lVert \prod\limits_{i=1}^{k-1}g^i\left(g^k_n-g^k\right)\right\rVert_\frac{p}{k}\\
        &\leq\left\lVert g^k_n\right\rVert_p\left\lVert \prod\limits_{i=1}^{k-1}g^i_n-\prod\limits_{i=1}^{k-1}g^i\right\rVert_\frac{p}{k-1}+\left\lVert \prod\limits_{i=1}^{k-1}g^i\right\rVert_\frac{p}{k-1}\left\lVert g^k_n-g^k\right\rVert_p\\
        &<M\frac{\delta}{2M}+M\frac{\delta}{2M}=\delta
        \end{split}
        \end{equation*}
        Thus $\prod\limits_{i=1}^kg^i_n\overset{L^\frac{p}{k}}{\longrightarrow}\prod\limits_{i=1}^kg^i$ as needed.
\end{proof}
\begin{corollary}
\label{adjconv}
    Let $\left\{g_n\right\}_{n=1}^\infty\subset W^{1,n-1}$ with $g_n\overset{W^{1,n-1}}{\longrightarrow}g$ then $\adj\nabla g_n\overset{L^1}{\longrightarrow}\adj\nabla g$.
\end{corollary}
\begin{proof}
    It follows directly from the claim and the definitions of adjoint matrix and Sobolev norm.
\end{proof}
\begin{theorem}[$VMO$ Change of Variables]
    Given $f\in VMO_\varphi\left(\Omega\right)\cap W^{1,n-1}\left(\Omega\right)$ with $\partial_nf\in L^\infty\left(\Omega\right)$, and $g$ is a smooth function with compact and connected support on $\mathbb{R}^n\backslash\overline{\eim\left(\varphi\right)}$ then:
    \begin{equation}
\int\limits_\Omega g\left(f\left(x\right)\right)\det\nabla f\left(x\right)dx =\deg\left(f,\Omega,p\right)\int\limits_{\mathbb{R}^n}g\left(z\right)dz\tag{\ref{change of var}}
    \end{equation}
    For any $p\in\supp g$.
\end{theorem}
\begin{remark}
    The above theorem applies to $\tilde{F}$ as in Claim \ref{codim}.
\end{remark}
\begin{proof}
    By the assumptions on $\supp g$, and Theorem \ref{bound2}, there is some $M$ neighborhood of $\Gamma$ on $\Omega$, and there are some $c>0$, and $d>0$ such that $f_i=f_T*\eta_\frac{c}{i}$ is well defined for all $i\in\mathbb{N}$, and has $\fint\limits_B\left\lvert f_i\left(y\right)-p\right\rvert dy>d$ for all $p\in\supp g$, $B\in D_M$, and $i\in\mathbb{N}$. Thus by Claim \ref{propvmodeg}, and Lemma \ref{mollify}, for any $p\in\supp g$ we have:
    \begin{equation}
    \label{lim1}
    \lim\limits_{i\rightarrow\infty}\deg\left(f_i,\Omega,p\right)=\deg\left(f,\Omega,p\right).
    \end{equation}
    Notice that $f_T\in W^{1,n-1}\left(\Omega\cup V\right)$, thus $\nabla f_i=\nabla f*\eta_\frac{c}{i}=\left(\nabla f\right)_i$ almost everywhere on $\Omega$. This means that $f_i\rightarrow f$ on $W^{1,n-1}\left(\Omega\right)$, thus both $f_i\rightarrow f$ and $\nabla f_i\rightarrow\nabla f$ point-wise a.e.\ on $\Omega$ (after passing to subsequence), giving $g\circ f_i \det\nabla f_i\rightarrow g\circ f\det\nabla f$ point-wise a.e.\ on $\Omega$.\\
    Recall that $\det\nabla f=\sum\limits_{j=1}^n\left(\adj\nabla f\right)^n_j(\nabla f)_n^j$. Thus, by the triangle inequality, and the assumption that $\partial_nf\in L^\infty$ we have that:
    \begin{equation*}
    \begin{split}
        \left\lvert g\circ f_i \det\nabla f_i\right\rvert&\leq\left\lvert g\circ f_i \right\rvert\sum\limits_{j=1}^n\left\lvert\left(\adj\nabla f_i\right)^n_j(\nabla f_i)_n^j\right\rvert\\
        &=\left\lvert g\circ f_i \right\rvert\sum\limits_{j=1}^n\left\lvert\left(\adj\nabla f_i\right)^n_j\right\rvert\left\lvert\left((\nabla f)_n^j\right)_i\right\rvert\\
        &\leq\left\lVert g\right\rVert_\infty\sum\limits_{j=1}^n\left\lvert\left(\adj\nabla f_i\right)^n_j\right\rvert\left\lVert\left(\left(\nabla f\right)_n^j\right)_i\right\rVert_\infty\\
        &\leq\left\lVert g\right\rVert_\infty\sum\limits_{j=1}^n\left\lvert\left(\adj\nabla f_i\right)^n_j\right\rvert\left\lVert\left(\left(\nabla f\right)_n^j\right)\right\rVert_\infty
    \end{split}
    \end{equation*}
    Where the last step follow from Young's convolution inequality:
    $$\left\lVert\left(\left(\nabla f\right)_n^j\right)_i\right\rVert_\infty=\left\lVert\left(\nabla f\right)_n^j*\eta_\frac{c}{i}\right\rVert_\infty\leq\left\lVert\left(\nabla f\right)_n^j\right\rVert_\infty\left\lVert\eta_\frac{c}{i}\right\rVert_1=\left\lVert\left(\nabla f\right)_n^j\right\rVert_\infty$$
    But by Corollary \ref{adjconv} we have $$\left\lVert g\right\rVert_\infty\sum\limits_{j=1}^n\left\lvert\left(\adj\nabla f_i\right)^n_j\right\rvert\left\lVert\left((\nabla f)_n^j\right)\right\rVert_\infty\overset{L^1}{\longrightarrow}\left\lVert g\right\rVert_\infty\sum\limits_{j=1}^n\left\lvert\left(\adj\nabla f\right)^n_j\right\rvert\left\lVert\left((\nabla f)_n^j\right)\right\rVert_\infty$$
    Thus by Pratt's theorem we must have:
    \begin{equation}
    \label{lim2}
        g\circ f_i\det\nabla f_i\overset{L^1}{\rightarrow}g\circ f\det\nabla f
    \end{equation}
    Thus by Equation \eqref{lim1}, Equation \eqref{lim2}, and Theorem \ref{changeofvar} we have:
    \begin{equation*}
        \begin{split}
            \deg\left(f,\Omega,p\right)\int\limits_{\mathbb{R}^n} g\left(z\right)dz&=\lim\limits_{i\rightarrow\infty}\deg\left(f_i,\Omega,p\right)\int\limits_{\mathbb{R}^n} g\left(z\right)dz\\
            &=\lim\limits_{i\rightarrow\infty}\int\limits_\Omega g\left(f_i\left(x\right)\right)\det\nabla f_i\left(x\right)dx\\
            &=\int\limits_\Omega g\left(f\left(x\right)\right)\det\nabla f\left(x\right)dx 
        \end{split}
    \end{equation*}
    As needed.
\end{proof}
\printbibliography
\end{document}